\documentclass{amsart}
\usepackage{amscd}
\usepackage{amssymb}
\input xypic
\usepackage[colorinlistoftodos]{todonotes}
\xyoption{all}
\newdir{ >}{{}*!/-5pt/@{>}}

\newcommand{\en}{\operatorname{End}}

\newcommand{\hp}{\mathcal{H}}
\renewcommand{\r}{\rightarrow}
\newcommand{\kh}{\operatorname{KH}}
\newcommand{\alg}{\operatorname{Alg}}

\newcommand{\coli}{\operatorname{colim}}

\def\id{\operatorname{id}}

\def\cC{\mathcal C}
\def\cT{\mathcal T}

\def\cE{\mathcal E}

\def\cU{\mathcal U}

\def\cG{\mathcal G}

\def\cM{\mathcal M}

\def\cS{\mathcal S}

\def\ev{\mathrm{ev}}

\newcommand{\C}{\mathbb{C}}

\newcommand{\Z}{\mathbb{Z}}
\newcommand{\N}{\mathbb{N}}
\newcommand{\kk}{{\operatorname{kk}}}

\numberwithin{equation}{section}
\theoremstyle{plain}
\newtheorem{thm}[equation]{Theorem}
\newtheorem{cor}[equation]{Corollary}

\newtheorem{prop}[equation]{Proposition}

\newcommand{\comment}[1]{}  

\theoremstyle{definition}
\newtheorem{defn}[equation]{Definition}
\theoremstyle{remark}
\newtheorem{rem}[equation]{Remark}
\newtheorem{ex}[equation]{Example}

\begin{document}

\bibliographystyle{plain}

\title{Algebraic quantum kk-theory}
\author{Eugenia Ellis}
\email{eellis@fing.edu.uy}
\address{IMERL, Facultad de Ingenier\'\i a, Universidad de la Rep\'ublica \\
Julio Herrera y Reissig 565, 11.300, Montevideo, Uruguay}
\thanks{The author was partially supported by ANII, CSIC and PEDECIBA}

\begin{abstract}
Let $\mathcal{G}$ be an algebraic quantum group. We introduce an equivariant algebraic $kk$-theory for $\mathcal{G}$-module algebras.
We study an adjointness theorem related with smash product and trivial action. We also discuss a duality property.
\end{abstract}

\maketitle

\section{Introduction}

Algebraic $kk$-theory is a bivariant $K$-theory on the category of
$\ell$-algebras where $\ell$ is a commutative ring with unit, \cite{CT}.
For each pair $(A,B)$ of $\ell$-algebras a group $\kk(A,B)$ is defined.
A category $\mathfrak{KK}$ is obtained whose objects are
$\ell$-algebras and where the morphisms from $A$ to $B$ are the elements of the group $\kk(A,B)$. The category $\mathfrak{KK}$ is triangulated and there is
a canonical functor $j: \alg_{\ell}\r \mathfrak{KK}$ with universal
properties. These properties are algebraic homotopy invariance, matrix
invariance and excision. The main result from \cite{CT} says that if $A$ is an $\ell$-algebra then
$$
\kk(\ell, A) \simeq \kh (A)
$$
here $\kh$ is Weibel's homotopy K-theory defined in \cite{kh}.

In \cite{eakk} an equivariant algebraic kk-theory is introduced. This is a  bivariant K-theory on the category of $G$-algebras (i.e. $\ell$-algebras with an action of $G$) where $G$ is a countable group. For each pair $(A,B)$ of $G$-algebras  a group  $\kk^{G}(A,B)$ is defined.   A category $\mathfrak{KK}^{G}$ is obtained whose objects are $G$-algebras and where the morphisms from $A$ to $B$ are the elements of the group $\kk^{G}(A,B)$. The category $\mathfrak{KK}^{G}$ is triangulated and there is a canonical functor $j: G$-$\alg \r \mathfrak{KK}^{G}$ with universal properties. These properties are algebraic homotopy invariance, equivariant matrix invariance and excision.

Let $\ell=\mathbb{C}$ and $\mathcal{G}$ be an algebraic quantum group in the sense of Van Daele (\cite{vd}, \cite{vd2}, \cite{ahma}). An algebraic quantum group $(\mathcal{G},\Delta, \varphi)$ is a regular multiplier Hopf algebra with invariants. That means, $\mathcal{G}$ is an associative algebra with non-degenerate product, $\Delta: \mathcal{G}\rightarrow M( \mathcal{G}\otimes \mathcal{G})$ is a coproduct (here $M(\mathcal{G}\otimes \mathcal{G})$ denotes the multiplier algebra of $\mathcal{G}\otimes \mathcal{G}$) such that there exists an invertible antipode $S:\mathcal{G}\rightarrow \mathcal{G}$ and $\varphi$ is a left invariant functional.  We define a bivariant K-theory on the category of $\mathcal{G}$-module algebras. For each pair $(A,B)$ of $\mathcal{G}$-module algebras we define a group $\operatorname{kk}^{\mathcal{G}}(A,B)$ and consider the category $\mathfrak{KK}^{\mathcal{G}}$ whose objects are the $\mathcal{G}$-module algebras and the morphisms from $A$ to $B$ are the elements of $\operatorname{kk}^{\mathcal{G}}(A,B)$.
Consider the functor $j^{\mathcal{G}}:\mathcal{G}$-$\operatorname{Alg} \rightarrow \mathfrak{KK}^{\mathcal{G}}$ which at the level of objects is the identity and at the level of morphisms sends $f:A\rightarrow B$ to its class $[f]$ in $\operatorname{kk}^{\mathcal{G}}(A,B)$. The category $\mathfrak{KK}^{\mathcal{G}}$ is triangulated and $j^{\mathcal{G}}$ is an excisive, homotopy invariant and weak $\mathcal{G}$-stable functor. Moreover, it is the universal functor for these properties.

The Green-Julg Theorem in Kasparov KK-theory (see \cite{Mkk}) states that if $G$ is a compact group and $B$ is a $G$-$C^*$-algebra, then there exists an isomorphism
$$
\mu : \operatorname{KK}^{G}(\C,B)\r\operatorname{KK}(\C, B\rtimes G).
$$
In the context of E-theory there is another version of Green-Julg Theorem (see \cite{ght}) which states that if $G$ is a compact group then  there is an isomorphism $$\mu: E_{G}(\C,B)\r E(\C,B\rtimes G).$$
Green-Julg Theorem in periodic cyclic homology is also studied in \cite{AK}, \cite{cv2} and \cite{cv1}.
In algebraic kk-theory there is a version of this theorem, see \cite{eakk}. Let $G$ be a finite group such that $\frac{1}{\mid G \mid} \in \ell$, $B$ an $\ell$-algebra and $A$ a $G$-algebra then there is an isomorphism
$$
\psi :\kk^{G}(B^{\tau},A)\r \kk (B,A\rtimes G)
$$
here $B^{\tau}$ denotes $B$ with the trivial action of $G$. In particular, if $B=\ell$ then
$$
\kk^{G}(\ell,A)\simeq \kk(\ell, A\rtimes G)\simeq \kh(A\rtimes G).
$$
The main theorem of this paper is a version of Green-Julg Theorem for $\mathfrak{KK}^{\hp}$ when $\hp$ is a semisimple Hopf algebra. Let $A$ be an algebra and $B$ an $\hp$-module algebra then there exists an isomorphism
$$
\psi :\kk^{\hp}(A^{\tau},B)\r \kk (A,B\#\hp)
$$
where $B\#\hp$ denotes the smash product.
In particular, if we take $A=\mathbb{C}$ we obtain that $\kk^{\hp}(\mathbb{C},B)\simeq \kh(\hp\#B)$.

The paper is organized as follows. In Section \ref{aqgsec} we recall the definition of an algebraic quantum group $\mathcal{G}$ and some of its properties. 
In Section \ref{akkth} we define algebraic $kk$-theory for $\mathcal{G}$-module algebras.
Let $A$ be a $\mathcal{G}$-module algebra.
In Section \ref{stability}  we present the stabilization property with respect to a triple $(V,W,\varphi)$ where $V$ and $W$ are left $\mathcal{G}$-modules and $\varphi: V\otimes W \rightarrow A$ is a $\mathcal{G}$-equivariant bilinear form. We introduce a strong $\mathcal{G}$-stabilization which imply equivariant Morita equivalence. Also we state a weak $\mathcal{G}$-stabilization which generalize the definition of $G$-stabilization with $G$ a group.     
In Section \ref{ekk} we define the algebraic quantum $kk$-theory and its properties. 
In Section \ref{gjekk} we prove Green-Julg Theorem when $\mathcal{G}=\mathcal{H}$ a semisimple Hopf algebra.
In Section \ref{bsd} we state a Baaj-Skandalis duality in this context. 

We would like to thank to the referee for his/her comments about the first manuscript of this paper.

\section{Algebraic quantum groups} \label{aqgsec}
In this section we recall the definition of algebraic quantum groups introduced in \cite{vd2}. Let $\mathcal{G}$ be an associative algebra over $\mathbb{C}$ with a non-degenerate product. This means that, if $hk=0$ for all $h\in \mathcal{G}$ then $k=0$ and if $hk=0$ for all $k\in \mathcal{G}$ then $h=0$.
A {\emph{multiplier}} of $\mathcal{G}$ is a pair $(\rho_{1},\rho_{2})$ such that $\rho_{i}:\mathcal{G}\rightarrow \mathcal{G}$ is a linear map ($i=1,2$) and
$$\rho_{1}(hk)=\rho_{1}(h)k\qquad \rho_{2}(hk)=h \rho_{2}(k)   \qquad \rho_{2}(h)k=h \rho_{1}(k) \qquad \qquad \forall h,k\in \mathcal{G}$$
We consider the associative natural product in the set of multipliers of $\mathcal{G}$,
$$(\rho_{1},\rho_{2})(\tilde{\rho}_{1},\tilde{\rho}_{2})=(\rho_{1}\tilde{\rho}_{1},\tilde{\rho}_{2}\rho_{2}).$$
Denote by $M(\mathcal{G})$ to the multiplier algebra of $\mathcal{G}$.
There is a natural embedding $\iota_{\mathcal{G}}: \mathcal{G}\rightarrow M(\mathcal{G})$,
$$
h\mapsto (L_{h},R_{h}) \qquad  L_{h}(k)=hk \qquad R_{h}(k)=kh
$$
which is a homomorphism.  Moreover $\rho h \in \mathcal{G} $ and $h\rho \in \mathcal{G}$ for all $h\in \mathcal{G}$ and $\rho \in M(\mathcal{G})$,
$$
\rho h =(L_{\rho_{1}(h)}, R_{\rho_{1}(h)}) \qquad  h\rho  =(L_{\rho_{2}(h)}, R_{\rho_{2}(h)})
$$
We write $\rho h =\rho_{1}(h)$ and  $ h\rho =\rho_{2}(h)$.
A {\emph{comultiplication}} on $\mathcal{G}$ is a homomorphism $\Delta: \mathcal{G}\rightarrow M(\mathcal{G}\otimes \mathcal{G})$ such that
$$
\Delta(h)(1\otimes k) \in \mathcal{G}\otimes \mathcal{G} \qquad (h\otimes 1)\Delta(k)\in \mathcal{G}\otimes \mathcal{G} \quad \forall h,k\in \mathcal{G}
$$
and the coassociativity property is satisfied
$$
(h\otimes 1 \otimes 1)(\Delta \otimes \id_{\mathcal{G}})(\Delta(k)(1\otimes r))=(\id_{\mathcal{G}}\otimes \Delta)((h\otimes 1)\Delta(k))(1\otimes 1\otimes r) \qquad \forall h,k,r \in \mathcal{G}
$$
A pair $(\mathcal{G},\Delta)$ is a {\emph{multiplier Hopf algebra}} if $\mathcal{G}$ is an algebra over $\mathbb{C}$ with a non-degenerate product, $\Delta: \mathcal{G} \rightarrow M(\mathcal{G}\otimes \mathcal{G})$ is a comultiplication on $\mathcal{G}$ and the following maps are bijective
$$
T_{i}:\mathcal{G}\otimes \mathcal{G} \rightarrow \mathcal{G} \otimes \mathcal{G} \qquad T_{1}(h\otimes k)=\Delta(h)(1\otimes k) \qquad T_{2}(h\otimes k)=(h\otimes 1) \Delta(k)
$$
If $(\mathcal{G},\Delta)$ is a multiplier Hopf algebra there is a unique homomorphism $\epsilon: \mathcal{G} \to \C$, called {\emph{counit}}, such that
$$
(\epsilon \otimes \id_{\mathcal{G}})(\Delta(h)(1\otimes k))=hk \qquad (\id_{\mathcal{G}}\otimes \epsilon)((h\otimes 1)\Delta(k))=hk \qquad \forall h,k\in \mathcal{G}.
$$
There is also a unique anti-homomorphism $S:\mathcal{G} \to M(\mathcal{G})$, called {\emph{antipode}}, such that
$$
m(S\otimes \id_{\mathcal{G}})(\Delta(h)(1\otimes k))=\epsilon(h)k
\qquad m(\id_{\mathcal{G}}\otimes S)((h\otimes
1)\Delta(k))=\epsilon(k)h   \qquad \forall h,k \in \mathcal{G}
$$
here $m$ is the multiplication map.
If $S(\mathcal{G})\subset \mathcal{G}$ and $S$ is invertible, $(\mathcal{G},\Delta)$ is called \emph{regular}. In this case we obtain that
$$
(1\otimes k) \Delta(h)\in \mathcal{G}\otimes \mathcal{G} \qquad \Delta(k)(h\otimes 1)\in \mathcal{G}\otimes \mathcal{G} \quad \forall h,k\in \mathcal{G}
$$
A {\emph{right invariant functional}} on $\mathcal{G}$ is a non-zero linear map $\psi:\mathcal{G} \to \C$ such that $$(\psi\otimes \id_{\mathcal{G}})\Delta(h)=\psi(h)1$$
Here $(\psi\otimes \id_{\mathcal{G}})\Delta(h)$ denotes the element $\rho \in M(\mathcal{G})$ such that
$$
\rho k = (\psi\otimes \id_{\mathcal{G}})(\Delta(h)(1\otimes k)) \qquad k\rho = (\psi\otimes \id_{\mathcal{G}})((1\otimes k)\Delta(h))
$$
Similarly, a {\emph{left invariant functional}} on  $\mathcal{G}$ is a non-zero linear map $\varphi:\mathcal{G} \to \C$ such that $$(\id_{\mathcal{G}}\otimes \varphi)\Delta(h)=\varphi(h)1$$
Invariant functionals do not always exist. If $\varphi$ is a left invariant functional on $\mathcal{G}$ then it is unique up to scalar multiplication and $\psi=\varphi \circ S$ is a right invariant functional.
Following \cite{ahma}, an {\emph{algebraic quantum group}} $(\mathcal{G}, \Delta)$ is a regular multiplier Hopf algebra with invariant funtionals on $\mathcal{G}$.

We will use the fact that algebraic quantum groups have local units, see  \cite[Proposition 2.6]{ahma}. That means given elements $\{h_{1},\ldots h_{n}\}$ in $\mathcal{G}$, there exists an element $e\in\mathcal{G}$ such that $h_{i}e=eh_{i}=h_{i}$ for all $i=1,\ldots n$.

Let $(\mathcal{G}, \Delta)$ be an algebraic quantum group. There exists an algebraic quantum group $(\hat{\mathcal{G}},\hat{\Delta})$ which is called the  {\emph{dual}} of  $(\mathcal{G}, \Delta)$. The elements of $\hat{\mathcal{G}}$ are the linear functionals of the form
$\varphi(h\cdot )$
$$
\hat{\mathcal{G}}=\{ \xi_{h}:\mathcal{G} \rightarrow \C: \xi_{h}(x)=\varphi(hx)\}
$$
The elements of $\hat{\mathcal{G}}$ can also be written as $\varphi(\cdot h)$, $\psi(h\cdot )$, $\psi(\cdot h)$. Write
$$
\begin{array}{lll}
\hat{\mathcal{G}}&=& \{ \eta_{h}:\mathcal{G} \rightarrow \C: \eta_{h}(x)=\psi(xh)\} \\
\\
&=&  \{ \hat{\eta}_{h}:\mathcal{G} \rightarrow \C: \hat{\eta}_{h}(x)=\psi(hx)\} \\
\\
&=& \{ \hat{\xi}_{h}:\mathcal{G} \rightarrow \C: \hat{\xi}_{h}(x)=\varphi(xh)\}
\end{array}
$$
The product on $\hat{\mathcal{G}}$ is defined as follows
$$
(\xi_{h}\cdot \xi_{k})(x)=(\varphi\otimes \varphi)(\Delta(x)(h\otimes k))
$$
The coproduct $\hat{\Delta}:\hat{\mathcal{G}}\to M(\hat{\mathcal{G}} \otimes \hat{\mathcal{G}})$ is defined by defining the elements $\hat{\Delta}(\xi_{1})(1\otimes \xi_{2})$ and $(\xi_{1}\otimes 1)\hat{\Delta}(\xi_{2})$ in $\hat{\mathcal{G}} \otimes \hat{\mathcal{G}}$ as follows
$$
((\xi_{1}\otimes 1)\hat{\Delta}(\xi_{2}))(h\otimes k)=(\xi_{1}\otimes \xi_{2})(\Delta(h)(1\otimes k))
$$
$$
(\hat{\Delta}(\xi_{1})(1 \otimes \xi_{2}))(h\otimes k)= (\xi_{1}\otimes \xi_{2})((h\otimes 1)(\Delta(k))
$$
The dual of $(\hat{\mathcal{G}},\hat{\Delta})$ is isomorphic to  $(\mathcal{G}, \Delta)$, see \cite[Theorem 4.12]{vd2}.

Let $G$ be a group then both $\mathcal{G}=\C G$ and its dual are algebraic quantum groups, see Example \ref{aqgg} and Example \ref{aqdg}. Another example of algebraic quantum groups are the finite-dimensional Hopf algebras as such algebras have invariants, see \cite[Prop 5.1]{vd2}.

Let $(\mathcal{G}, \Delta)$ be a regular multiplier algebra.
A {\emph{left $\mathcal{G}$-module}} $M$ is a vector space over $\C$ with a linear map $\mathcal{G}\otimes M\to M$, $h\otimes m \mapsto h\cdot m$ such that $(hk)\cdot m = h\cdot(k\cdot m)$ for all $h, k \in \mathcal{G}$ and $m\in M$. We say $M$ is {\emph{unital}} if $\mathcal{G}\cdot M=M$.
This condition coincides with the condition $1\cdot m =m$ when $\mathcal{G}$ has an identity.

Let $A$ be an algebra over $\C$ with or without identity but with
non-degenerate product. Assume that $A$ is a  {\emph{firm algebra}}, i.e. the
map $$ A\otimes_{A}A \rightarrow A \qquad r\otimes s \mapsto rs
$$ is an isomorphism. We say thay $A$ is a {\emph left $\mathcal{G}$-module
algebra} if $A$ is a unital left $\mathcal{G}$-module such that
\begin{equation}\label{mulal}
h\cdot(ab)=m(\Delta(h)(a\otimes b)) \qquad \forall a,b\in A, \ h\in \mathcal{G}
\end{equation}
It is possible to use a generalized Sweedler notation in the following way. The element $\Delta(h)$ is not in $\mathcal{G}\otimes \mathcal{G}$ but $\Delta(h)(1\otimes k)$ lies in $\mathcal{G}\otimes \mathcal{G}$ for all $k\in \mathcal{G}$.
If we fix a finite set of elements in $\mathcal{G}$, $\{k_{1},\ldots k_{n}\}$, there exists an element $e\in \mathcal{G}$ such that $ek_{i}=k_{i}=k_{i}e$ for all $i=1,\ldots, n$.  If we denote
$\sum h_{(1)}\otimes h_{(2)}$ for $\Delta(h)(1\otimes e)$ we can write
$$
\sum h_{(1)}\otimes h_{(2)}k_{i}=\Delta(h)(1\otimes k_{i})\quad \sum h_{(1)}\otimes k_{i} h_{(2)}=(1\otimes k_{i})\Delta(h)
$$
$$
\sum k_{i} h_{(1)}\otimes h_{(2)}=(k_{i}\otimes 1)\Delta(h)\quad \sum h_{(1)}k_{i}\otimes  h_{(2)}=\Delta(h)(k_{i}\otimes 1)
$$
Using this notation we rewrite \eqref{mulal} as
$$
h\cdot(ab)=\sum (h_{(1)}\cdot a)(h_{(2)}\cdot b)
$$
Because $A$ is unital there is $e\in \mathcal{G}$ such that $(a\otimes b)=(a\otimes e\cdot b).$

\begin{ex}\label{aqgg}
Let $G$ be a group and denote $e=1_{G}$. Take $\mathcal{G}=\C G$ with the usual Hopf algebra structure. The maps involved are the $\C$-linear extension of the following maps
$$
\begin{array}{c}
\mathcal{G}=\C G=\{\displaystyle \sum_{g\in G}a_{g}g : \mbox{$a_{g}\in \C$ $a_{g}\neq 0$ for a finite amount of $g$}\}
\\
\\
\begin{array}{ll}
\Delta: G\rightarrow G \otimes G  & \Delta(g)=g\otimes g\\
\\
T_{1}:G\otimes G \rightarrow G\otimes G & T_{1}(h\otimes k)=h\otimes hk\\
\\
T_{2}:G\otimes G \rightarrow G\otimes G & T_{2}(h\otimes k)=hk\otimes k\\
\\
\epsilon: G\rightarrow \C & \epsilon(g)=1\\
\\
S:G\rightarrow G & S(g)=g^{-1}\\
\\
\varphi=\psi=\chi_{e}: G\rightarrow \C & \chi_{e}(h)=\left\{ \begin{array}{ll} 1 & e=h\\
0 & e\neq h
\end{array}\right.\\
\\
\hat{G}=\{\chi_{h^{-1}}:G \rightarrow \C:  \chi_{h^{-1}}(g)=\left\{ \begin{array}{ll} 1 & h^{-1}=g\\
0 &  h^{-1}\neq g
\end{array}\right.\}\\
\end{array}\\
\\
\hat{\mathcal{G}}=\C\hat{G}=\{\displaystyle \sum_{g\in G}a_{g}\chi_{g} : \mbox{$a_{g}\in \C$ $a_{g}\neq 0$ for a finite amount of $g$}\}
\end{array}
$$
\end{ex}
\begin{ex}\label{aqdg} In the previous example we see the usual Hopf algebra structure of $\C G$. The dual of it can be not a Hopf algebra, for example when $G$ is not finite. Nevertheless it is an algebraic quantum group with the following structure
$$
\begin{array}{c}
\mathcal{G}=\C \hat{G}=\{\displaystyle \sum_{g\in G}a_{g}\chi_{g} : \mbox{$a_{g}\in \C$ $a_{g}\neq 0$ for a finite amount of $g$}\} \\ \chi_{g}\chi_{h}=  \left\{
\begin{array}{ll}
\chi_{g} & g=h\\
0 & g\neq h
\end{array}\right.
\\
\\
\begin{array}{ll}
\Delta: \mathcal{G}\rightarrow M(\mathcal{G} \otimes \mathcal{G})  & \Delta(\chi_g)=\displaystyle \sum_{t\in G} \chi_{gt^{-1}}\otimes \chi_t\\
\\
\end{array}
\end{array}
$$
This sum is not necessarily finite but
$$\Delta(\chi_{h})(1\otimes \chi_{k})=\chi_{hk^{-1}}\otimes \chi_{k}\in \mathcal{G}\otimes\mathcal{G} \quad
(\chi_{h}\otimes 1)\Delta(\chi_{k})=\chi_{h}\otimes \chi_{h^{-1}k}\in\mathcal{G}\otimes\mathcal{G} $$
The counit and the antipode are
$$
\epsilon:\mathcal{G}\rightarrow \C \quad \epsilon(\chi_{g})=\left\{ \begin{array}{cc} 1 & g=e \\
0 & g\neq e
\end{array}\right. \qquad S:\mathcal{G}\rightarrow\mathcal{G} \quad S(\chi_{g})=\chi_{g^{-1}}
$$
The integral is
$$
\varphi=\psi:\mathcal{G}\rightarrow \C \qquad \varphi(\chi_{h})=\psi(\chi_{h})=1
$$
\end{ex}

\section{Algebraic $kk$-theory}\label{akkth}

Let $\mathcal{G}$ be an algebraic quantum group.
We define an algebraic $kk$-theory for the category of $\mathcal{G}$-module algebras. We write $\cC$ to refer to this category.
We adapt some results from \cite{CT} and \cite[Section 2]{eakk} to our setting.

Let $\mathcal{U}(\cC)$ be the category of $\mathcal{G}$-modules, i.e. we forget the multiplication in $\cC$ keeping the $\mathcal{G}$-module structure.
Let $F:\cC \r \mbox{$\mathcal{U}(\cC)$}$ be the forgetful functor.
Let $M$ be an object in $\cU(\cC)$. Take
$$
\tilde{T}(M)=\bigoplus_{n\geq 1} M^{\otimes^{n}}\qquad \qquad
M^{\otimes^{n}}=\underset{\mbox{$n$-times}}{\underbrace{M\otimes
\ldots \otimes M}}
$$
with the usual structure of  algebra. Consider in $M^{\otimes^{n}}$ the following action,
$$
\displaystyle{
h\cdot(m^{1}\otimes m^{2} \otimes \ldots \otimes m^{n})=
\sum h_{(1)}\cdot m^{1}\otimes h_{(2)}\cdot m^{2} \otimes \ldots \otimes h_{(n)}\cdot m^{n}
}
$$
which gives to $\tilde{T}(M)$ an $\mathcal{G}$-module algebra structure.
Both constructions are functorial hence we consider the functor
 $\tilde{T}:\cU(\cC)\r \cC$.
Put
$$
T:=\tilde{T}\circ F:\cC \r \cC.
$$
The functor $\tilde{T}$ is the left adjoint of $F$ and the counit
of the adjuntion is
$$\eta_{A}: T(A)\r A \qquad \eta_{A}(a_{1}\otimes \ldots \otimes
a_{n})=a_{1}\ldots a_{n}$$ and it is surjective (see
\cite[IV.3 Thm 1]{Mac}). Put $J(A)=\ker \eta_{A}$ and note this construction is functorial.
The universal extension of $A$ is
$$
J(A)\xrightarrow{\imath_{A}} T(A)\xrightarrow{\eta_{A}} A.
$$
Let $A\xrightarrow{f} B \xrightarrow{g} C$ be a weakly split extension, (i.e. $F(g)$ has a section in $\cU(\cC)$).
There is a map $\xi: J(C)\r A$ called classifying map of $A\xrightarrow{f} B \xrightarrow{g} C$. This map is defined as follows.  Let $s$ be a section of $F(g)$ and consider the map $\hat{\xi}:T(C)\rightarrow B$, $\hat{\xi}=\eta_{B}\circ \tilde{T}(s)$.
The map $\xi: J(C)\r A$ is the restriction of $\hat{\xi}$ to $J(C)$.
This construction of $\xi$ depends on which $s$ is chosen
but $\xi$ is unique up to elementary homotopy.

Let $L$ be a ring and $A$ an object in $\cC$. Then the extension
\begin{equation}\label{jaltal}
J(A)\otimes_{\Z} L\r T(A)\otimes_{\Z}L\r A\otimes_{\Z}L
\end{equation}
is weakly split, and there is a choice for the classifying map $$\phi_{A,L}:
J(A\otimes_{\Z}L)\r J(A)\otimes_{\Z}L$$ of \eqref{jaltal}, which is
natural in both variables.

Let $A$ and $B$ be objects in $\cC$ and $\cM_{\infty}$ the ind-ring defined in \cite[Sec 4.1]{CT}. Let $S^{1}$ be the simplicial circle $\Delta^{1}/\partial\Delta^{1}$ and $\operatorname{sd}^{\bullet}S^{1}$ the pro-simplicial set of subdivision of $S^{1}$.
  Recall that $A^{(\operatorname{sd}^{n} S^{1},\star)}$ is the algebra of polynomial functions on the pointed simplicial set $(\operatorname{sd}^{n} S^{1},\star)$.
By \cite[Lemma 3.1.3]{CT} we obtain that $A\otimes \Z^{(\operatorname{sd}^{n} S^{1},\star)}\simeq A^{(\operatorname{sd}^{n} S^{1},\star)}$ as  algebras. If $A$ is an $\mathcal{G}$-module algebra we can consider $A^{(\operatorname{sd}^{n} S^{1},\star)}$ as an $\mathcal{G}$-module algebra taking the diagonal structure with the trivial structure in $\Z^{(\operatorname{sd}^{n} S^{1},\star)}$. Define inductively
$$A^{\cS^{1}}:=A^{(\operatorname{sd}^{\bullet} S^{1},\star)} \quad  B^{\cS^{n+1}}:=(B^{\cS^{n}})^{\cS^{1}}.$$
The map $\rho_{A}:J(A)\r A^{\cS^{1}}$ denotes the composition of the classifying map of the loop extension
\begin{equation}{\label{extlaz}}
\Omega A\r PA\xrightarrow{\ev_{1}}A \qquad \Omega A= \Z^{(S^{1},\star)}\otimes A \qquad  P A= \Z^{(\Delta^{1},\star)}\otimes A
\end{equation}
with the last vertex map $\Omega A \xrightarrow{h} A^{\cS^{1}}$.
Consider the algebraic homotopy classes in $\cC$ of the maps from $J^{n}A$ to $\cM_{\infty} B^{\cS^{n}}$. Denote it by
$$
E_{n}(A,B)_{\cC}:=[J^{n}(A),\cM_{\infty} B^{\cS^{n}}]_{\cC}.
$$
Consider the following morphism
$\imath_{n}:E_{n}(A,B)_{\cC}\r E_{n+1}(A,B)_{\cC}$ such that
$$\begin{array}{c}
J^{n}(A)\xrightarrow{f}\cM_{\infty} B^{\cS^{n}}\\
\downarrow
\\
 J^{n+1}(A)
\xrightarrow{J(f)}
J(\cM_{\infty}B^{\cS^{n}})\xrightarrow{\phi_{\cM_{\infty},B^{\cS^{n}}}}
\cM_{\infty}J(B^{\cS^{n}})\xrightarrow{\rho_{B^{\cS^{n}}}}\cM_{\infty}B^{\cS^{n+1}}
\end{array}
$$
Define
$$
\kk_{\cC}(A,B)=\coli_{n\in \N}E_{n}(A,B)_{\cC}.
$$

Let $\mathfrak{KK}_{\cC}$ be the category whose objects are the same objects of $\cC$ and the morphisms from $A$ to $B$ are the elements of $\kk_{\cC}(A,B)$. The composition rule is given by the product defined in \cite[Section 6.2]{CT}.
Denote by
$$
j_{\cC}:\cC \r \mathfrak{KK}_{\cC}
$$
the functor which at the level of objects is the identity and at
level of morphisms sends $f:A\r B$ to $[f]\in \kk_{\cC}(A,B)$. Consider the functor $\Omega:\mathfrak{KK}_{\cC} \r \mathfrak{KK}_{\cC}$ which sends an object
  $A$ of $\cC$  to the path object $\Omega A$. A diagram
$$
\Omega C \r A\r B\r C
$$
of morphisms in $\mathfrak{KK}_{\cC}$ is called a distinguished
triangle if it is isomorphic in $\mathfrak{KK}_{\cC}$ to the path sequence associated with a homomorphism $f:A'\to B'$ in $\cC$
\begin{equation}\label{psec}
\Omega B'\xrightarrow{j(\iota)} P_{f}\xrightarrow{j(\pi_{f})}
A'\xrightarrow{j(f)} B'
\end{equation}
Recall that \eqref{psec} is obtained from the loop extension of $B'$ by pulling it back to $A'$
$$
P_{f}:= PB' \times_{B}A'\quad\quad \xymatrix{
\Omega B' \ar[r]^-{\iota}\ar[d] & PB' \times_{B'}
A'\ar[r]^-{\pi_{f}}\ar[d] & A'\ar[d]^{f}& \operatorname{E}'\\
\Omega B'\ar[r]_-{\iota}& \ar[r] PB' \ar[r]_-{\ev_{1}}& B'& \operatorname{E}
}
$$
The category $\mathfrak{KK}_{\cC}$ is triangulated with the structure defined above, see \cite[Theorem 6.5.2]{CT}.

Let $\operatorname{E}: A\xrightarrow{f}B\xrightarrow{g} C$ be a weakly split
extension and let
$c_{\operatorname{E}}\in \kk_{\cC}(J(C),A)$ its classifying map.
As the natural map $\rho_{A}:J(A)\r \Omega A$ induces a
$\kk_{\cC}$-equivalence (see Lemma 6.3.10, \cite{CT}) we can consider $\partial_{\operatorname{E}}:=c_{\operatorname{E}}\circ \rho^{-1}_{C}$ in $\kk_{\cC}(\Omega C, A)$.

\begin{thm}\label{pukk}
The functor $j_{\cC}:\cC\r \mathfrak{KK}_{\cC}$ with the morphisms
 $\{\partial_{\operatorname{E}}:\operatorname{E}\in\cE\}$ is an excisive homology theory (in the sense of \cite[Section 6.6]{CT}), homotopy invariant and
$M_{\infty}$-stable (in the sense of \cite[Section 2.1, 2.2]{eakk}).  Moreover, it is universal with these properties. In other words,  if $\cT$
is a triangulated category and $G:\cC \r \cT$ together with
a class of morphisms $\{\overline{\partial}_{\operatorname{E}}:\operatorname{E}\in\cE\}$ is an
excisive, homotopy invariant and $M_{\infty}$-stable functor,
then there exists a unique triangle functor $\overline{G}:
 \mathfrak{KK}_{\cC}\r \cT$ such that the following diagram commutes
$$
\xymatrix{
\cC \ar[r]^{j_{\cC}}\ar[dr]_{G}& \mathfrak{KK}_{\cC}\ar@![d]^{\overline{G}}\\
&\cT
}
$$
\end{thm}
\proof See \cite[Thm 2.6.5]{eakk} and \cite[Thm 6.6.6]{CT}.

Let $\mathcal{X}$ be a infinity set. Consider
$$M_{\mathcal{X}}=\{ f:\mathcal{X}\times \mathcal{X}\rightarrow \C: \qquad \operatorname{supp}(f)<\infty \}$$
Consider the category $\mathfrak{KK}_{\mathcal{X}}$ with the same objects of  $\mathfrak{KK}_{\cC}$ and
$$
\mathfrak{KK}_{\mathcal{X}}(A,B):=\mathfrak{KK}_{\cC}(A,M_{\mathcal{X}}B)
$$

We take the following result

\begin{thm}\label{xpuk}\cite[Theorem 5.2.16]{RC}\cite[Theorem 9.3(2)]{Gar2}
The functor $j_{\cC}:\cC\r \mathfrak{KK}_{\mathcal{X}}$ with the morphisms
 $\{\partial_{\operatorname{E}}:\operatorname{E}\in\cE\}$ is an excisive homology theory, homotopy invariant and
$M_{\mathcal{X}}$-stable (in the sense of \cite[5.2.10]{RC}). It is universal for these properties.
\end{thm}

\section{Stability}\label{stability}

\subsection{Triples} Let $V$, $W$ be vector spaces over $\mathbb{C}$ with a bilinear form $\phi:V\otimes W\rightarrow
\mathbb{C}$. Define an algebra
$W\underset{\tiny\phi}{\otimes}V$ :
$$  W\underset{\tiny\phi}{\otimes}V =W{\otimes}V
\quad (w\otimes v)(\tilde{w}\otimes \tilde{v})= w\otimes \phi (v
\otimes \tilde{w}) \tilde{v}.$$ Let $(\mathcal{G},\Delta)$ be an
algebraic quantum group and $A$ be a $\mathcal{G}$-module algebra.
We say that $(V,W,\phi)$ is a {\emph{$(\mathcal{G},A)$-triple}} if
\begin{itemize}
\item $V$ is an unital left $\mathcal{G}$-module and a left $A$-module such that
$$g\cdot(a\cdot v)=\sum(g_{(1)}\cdot a)\cdot (g_{(2)}\cdot v) \qquad \qquad \forall g\in \mathcal{G} \quad a\in A \quad v\in V$$
\item $W$ is an unital left $\mathcal{G}$-module and a right $A$-module such that
$$g\cdot(w \cdot a)=\sum(g_{(1)}\cdot w)\cdot (g_{(2)}\cdot a) \qquad \qquad \forall g\in \mathcal{G} \quad a\in A \quad w\in W$$
\item $\phi: V\otimes W \rightarrow A$ is a $\mathcal{G}$-equivariant $A$-bimodule.
\end{itemize}
Let $(V,W,\phi)$ be a $(\mathcal{G},A)$-triple.
We obtain that $W\otimes_{A} V$ is an associative algebra with the following product
 $$(w_{1}\otimes v_{1})(w_{2}\otimes v_{2})= w_{1}\otimes \phi(v_{1}\otimes
w_{2})\cdot v_{2}=  w_{1}\cdot\phi(v_{1}\otimes w_{2}) \otimes v_{2}
 $$
 We denote it as $W\underset{\tiny\phi}{\otimes}V$.
 We consider the following action on $W\underset{\tiny\phi}{\otimes}V$
 $$
 g\cdot
(w\otimes v)=\sum g_{(1)}\cdot w \otimes g_{(2)}\cdot v
 $$
We check that  $W\underset{\tiny\phi}{\otimes}V$ is a $\mathcal{G}$-module algebra with this structure:
$$
\begin{array}{lll}
g\cdot((w_{1}\otimes v_{1})(w_{2}\otimes v_{2}))&=&g\cdot (w_{1}\otimes \phi(v_{1}\otimes w_{2})\cdot v_{2})\\
&=&\sum g_{(1)}\cdot w_{1}\otimes  g_{(2)}\cdot (\phi( v_{1}\otimes
w_{2})\cdot v_{2})\\
&=& \sum g_{(1)}\cdot w_{1}\otimes \phi(g_{(2)}\cdot v_{1}\otimes
g_{(3)}\cdot w_{2})\cdot (g_{(4)}\cdot v_{2})\\
 &=& \sum g_{(1)}\cdot
(w_{1}\otimes v_{1}) g_{(2)}\cdot (w_{2}\otimes  v_{2})
\end{array}
$$
If $A=\C$ we say that $(V,W,\phi)$ is a  $\mathcal{G}$-triple.
\begin{ex}
Recall that
$$
\mathcal{G}^{*}=\{ \eta: \mathcal{G} \rightarrow \C: \mbox{$\eta$ is
$\C$-linear}\}\qquad \hat{\mathcal{G}}=\{ \xi_{h}:\mathcal{G}
\rightarrow \C: \xi_{h}(x)=\varphi(hx)\}
$$
Consider the regular action on $\mathcal{G}$, $\operatorname{ev}:
{\mathcal{G}^{*}}\otimes \mathcal{G}\to \C$ and
$\operatorname{ev}:\hat{\mathcal{G}}\otimes \mathcal{G}\to \C$ the
evaluation maps and the following action of $\mathcal{G}$ on
${\mathcal{G}^{*}}$ and $\hat{\mathcal{G}}$
$$
(g\cdot \eta)(h)=\eta (S(g)h) \qquad g\cdot \xi_{k}=\xi_{k S(g)}
$$
Consider the following two particular $\mathcal{G}$-triples:
$({\mathcal{G}^{*}},\mathcal{G}, \operatorname{ev} )$ and
$(\hat{\mathcal{G}},\mathcal{G}, \operatorname{ev} )$. Define $\mathcal{G}\diamondsuit {\mathcal{G}^{*}}$ and  $\mathcal{G}\diamondsuit \hat{\mathcal{G}}$ as the following $\mathcal{G}$-module algebras
$$
\mathcal{G}\diamondsuit {\mathcal{G}^{*}} =
\mathcal{G}\underset{\tiny\operatorname{ev}}{\otimes}\mathcal{G}^{*}\qquad
\qquad \mathcal{G}\diamondsuit \hat{\mathcal{G}} =
\mathcal{G}\underset{\tiny\operatorname{ev}}{\otimes}\hat{\mathcal{G}}
$$
It is easy to see that $\mathcal{G}\diamondsuit \hat{\mathcal{G}}$
is a $\mathcal{G}$-module subalgebra of $\mathcal{G}\diamondsuit
{\mathcal{G}^{*}}$.
\end{ex}

\begin{ex}
Let $A$ be a $\mathcal{G}$-module algebra. We can consider the
triple $(A,A,m_{A})$ with the left and right regular action on $A$
and $m_{A}:A\otimes A\rightarrow A$ the multiplication on $A$.
\end{ex}
Given a $(\mathcal{G},A)$-triple  $(V,W,\phi)$,  we construct
$(V_{\mathcal{G}}, W_{\mathcal{G}}, \phi_{\mathcal{G}})$ a
$(\mathcal{G},A)$-triple associated to $(V,W,\phi)$ by the following
way
$$
V_{\mathcal{G}}=V\otimes \hat{\mathcal{G}} \qquad W_{\mathcal{G}}=\mathcal{G}\otimes W \qquad \phi_{\mathcal{G}}: V_{\mathcal{G}}\otimes W_{\mathcal{G}}\xrightarrow{\operatorname{ev}} V\otimes W \xrightarrow{\phi} A
$$
The left $\mathcal{G}$-module structures on $V_{\mathcal{G}}$ and
$W_{\mathcal{G}}$ are given by the diagonal action:
$$
g\cdot(v\otimes \xi_{k})= \sum g_{(1)}\cdot v\otimes
\xi_{kS(g_{(2)})} \qquad g\cdot(h\otimes w)=\sum g_{(1)}h\otimes
g_{(2)}\cdot w.
$$

\begin{prop}\label{motmg} Let $(V,W,\phi)$ be a $(\mathcal{G},A)$-triple then
$$
W_{\mathcal{G}}\underset{\tiny\phi_{\mathcal{G}}}{\otimes}V_{\mathcal{G}}\simeq
\mathcal{G}\diamondsuit \hat{\mathcal{G}} \otimes
W\underset{\tiny\phi}{\otimes}V
$$
as algebras.
\end{prop}

\subsection{The algebra $\mathcal{A}(M)$}
Let $M$ be a $\mathcal{G}$-module. Let consider $(V_{M}, W_{M}, \phi_{M})$ the following $\mathcal{G}$-triple associated to $M$
$$
V_{M}= M^{*}=\{\eta: M\to \C : \mbox{$\eta$ is $\C$-linear }\} \qquad W_{M}=M \qquad \phi_{M}=\operatorname{ev}:M^{*}\otimes M \to \C
$$
The $\mathcal{G}$-module structure on $M^{*}$ is
$$
(g\cdot \eta)(m)=\eta(S(g)\cdot m)
$$
We define the $\mathcal{G}$-module algebra $\mathcal{A}(M)$ associated to $M$ as
$$
\mathcal{A}(M):= W_{M}\underset{\tiny\phi_{M}}{\otimes}V_{M}=  M\underset{\tiny\operatorname{ev}}{\otimes}M^{*}
$$
\begin{prop}\label{misoa}
Let $f:M\rightarrow N$ an isomorphism of $\mathcal{G}$-modules, then
$f\otimes (f^{-1})^{*}:\mathcal{A}(M)\rightarrow \mathcal{A}(N) $ is
an isomorphism of $\mathcal{G}$-module algebras.
\end{prop}

We say that $M$ has a {\emph{good duality}} if there is a submodule
$\hat{M}$ of $M^{*}$ together a $\mathcal{G}$-equivariant module
isomorphism:
$$
\omega: M \rightarrow \hat{M}\subset M^{*}
$$
Define $\hat{\mathcal{A}}(M)$ as the following subalgebra of
$\mathcal{A}(M)$
$$
\hat{\mathcal{A}}(M)=
M\underset{\tiny\operatorname{ev}}{\otimes}\hat{M} \subset
M\underset{\tiny\operatorname{ev}}{\otimes}M^{*} =\mathcal{A}(M)
$$
\begin{ex}
The algebraic quantum group $\mathcal{G}$ is a $\mathcal{G}$-module
with good duality with the regular action. The following map is an
equivariant isomorphism
$$
\omega: \mathcal{G} \rightarrow \hat{\mathcal{G}} \qquad
\omega(g)=\xi_{S(g)}
$$
\end{ex}

Let $T=(V, W, \phi)$ be a $\mathcal{G}$-triple and $A$ be a $\mathcal{G}$-module algebra. We construct the following $(\mathcal{G},A)$-triple
$$
V_{T,A}=A\otimes V \quad W_{T,A}=W\otimes A \quad \phi_{T,A}: A\otimes V\otimes W\otimes A\xrightarrow{\phi}A\otimes A \xrightarrow{m_{A}} A
$$

\begin{prop}
 Let $T=(V, W, \phi)$ be a $\mathcal{G}$-triple and $A$ be a $\mathcal{G}$-module algebra. Then
$$
 W_{T,A}\underset{\tiny\phi_{T,A}}{\otimes}V_{T,A} \simeq  A\underset{\tiny A}{\otimes}A \otimes W\underset{\tiny\phi}{\otimes}V \simeq  A \otimes W\underset{\tiny\phi}{\otimes}V
$$
as algebras.
\end{prop}

\begin{prop}
Let $M$ be a left $\mathcal{G}$-module, then $\mathcal{G}\otimes M$ is a $\mathcal{G}$-module with the following action
$$
g\cdot (k\otimes m) = \sum g_{(1)}k\otimes g_{(2)}\cdot m
$$
If $M^{\tau}$ denotes the trivial module, then  $\mathcal{G}\otimes M^{\tau} \simeq \mathcal{G}\otimes M$.
\begin{proof}
$$
\begin{array}{c}
\xymatrix{
\mathcal{G}\otimes M^{\tau} \ar@/^0.3pc/[r]^{\alpha}  &\ar@/^0.3pc/[l]^{\beta} \mathcal{G}\otimes M
}
\\
\alpha(g\otimes m)=\sum g_{(1)}\otimes g_{(2)}\cdot m \qquad \beta(g\otimes m)=\sum g_{(1)}\otimes S(g_{(2)})\cdot m
\end{array}
$$
\end{proof}
\end{prop}

\begin{prop}\label{modmodt}
Let $M$ be a left $\mathcal{G}$-module, then $\hat{\mathcal{G}}\otimes M$ is a $\mathcal{G}$-module with the following action
$$
g\cdot (\xi_{k}\otimes m) = \sum\xi_{k S( g_{(1)})}\otimes g_{(2)}\cdot m
$$
If $M^{\tau}$ denotes the trivial module, then  $\hat{\mathcal{G}}\otimes M^{\tau} \simeq \hat{\mathcal{G}}\otimes M$.
\begin{proof}
$$
\begin{array}{c}
\xymatrix{
\hat{\mathcal{G}}\otimes M^{\tau} \ar@/^0.3pc/[r]^{\alpha}  &\ar@/^0.3pc/[l]^{\beta} \hat{\mathcal{G}}\otimes M
}
\\
\alpha(\xi_{g}\otimes m)=\sum \xi_{g_{(2)}}\otimes S^{-1}(g_{(1)})\cdot m \qquad \beta(\xi_{g}\otimes m)=\sum \xi_{g_{(2)}}\otimes g_{(1)}\cdot m
\end{array}
$$
\end{proof}
\end{prop}

Let $M$ be a $\mathcal{G}$-module. Observe that
$\mathcal{A}(\mathcal{G})\otimes \mathcal{A}(M)$ is not a
$\mathcal{G}$-module algebra with the diagonal structure. Despite of
this, we can redefine the action on $\mathcal{A}(\mathcal{G})\otimes
\mathcal{A}(\mathcal{M})$ in order to obtain a $\mathcal{G}$-module
algebra.

\begin{prop} Let $M$ a $\mathcal{G}$-module.  The algebra $\mathcal{A}(\mathcal{G})\otimes
\mathcal{A}(M)$ is a $\mathcal{G}$-module algebra with the usual
product
$$
(g\otimes \eta \otimes m \otimes f) (\tilde{g}\otimes \tilde{\eta}
\otimes \tilde{m} \otimes \tilde{f})= (g\otimes
\eta)(\tilde{g}\otimes \tilde{\eta})\otimes (m \otimes f) (\tilde{m}
\otimes \tilde{f})=\eta(\tilde{g})f(\tilde{m}) g\otimes \tilde{\eta}
\otimes m \otimes \tilde{f}$$
 and the following action:
$$
t\cdot (g\otimes \eta \otimes m \otimes f)= \sum t_{(1)}\cdot g\otimes t_{(4)}\cdot \eta \otimes t_{(2)}\cdot m \otimes t_{(3)}\cdot f
$$
\end{prop}

\begin{prop}
Let $M$ be a $\mathcal{G}$-module, then there is an injection
$$
 \mathcal{A}({\mathcal{G}}) \otimes \mathcal{A}(M) \to \mathcal{A}(\mathcal{G}\otimes M)
$$
of $\mathcal{G}$-modules algebras.
\begin{proof}
If $\eta \in {\mathcal{G}}^{*}$ and $f \in M^{*}$ then $\eta \otimes f \in (\mathcal{G}\otimes M)^{*}$. Then we have the following injection which is a morphism of $\mathcal{G}$-module algebras
$$
\mathcal{G}\diamondsuit {\mathcal{G}}^{*} \otimes M\underset{\tiny\operatorname{ev}}{\otimes}M^{*}\to  (\mathcal{G}\otimes M)\underset{\tiny\operatorname{ev}}{\otimes}(\mathcal{G} \otimes M)^{*} \quad g\otimes \eta \otimes m \otimes  f \mapsto  g \otimes m\otimes  \eta \otimes f
$$
\end{proof}
\end{prop}

We would like to have an isomorphism between $\hat{\mathcal{A}}(\mathcal{G})\otimes\mathcal{A}(M)$ and  $\hat{\mathcal{A}}(\mathcal{G})\otimes\mathcal{A}(M^{\tau})$.
In the following proposition we write an explicit isomorphism.

\begin{prop}\label{isowow}
 Let $M$ be a $\mathcal{G}$-module, then the following map is a $\mathcal{G}$-module algebra isomorphism
$$
\begin{array}{llll}
\delta: & \hat{\mathcal{A}}({\mathcal{G}})\otimes \mathcal{A}(M)& \to& \hat{\mathcal{A}}({\mathcal{G}})\otimes \mathcal{A}(M^\tau)\\ \\ & g\otimes \eta_{h} \otimes m \otimes  f &\mapsto &\sum g_{(1)}\otimes   \eta_{h_{(1)}} \otimes S(g_{(2)})\cdot  m\otimes h_{(2)}\cdot f
\end{array}
$$
\begin{proof}
Take $a= g\otimes \eta_{h} \otimes m \otimes  f$ and $b= \tilde{g}\otimes \eta_{\tilde{h}} \otimes \tilde{m} \otimes  \tilde{f}$ then
$$
ab=g\eta_{h}(\tilde{g})\otimes \eta_{\tilde{h}}\otimes m f(\tilde{m})\otimes \tilde{f} \qquad \delta{(ab)}=\sum \eta_{h}(\tilde{g})f(\tilde{m}) g_{(1)}\otimes \eta_{\tilde{h}_{(1)}}\otimes S(g_{(2)})\cdot m \otimes \tilde{h}_{(2)}\cdot \tilde{f}.
$$
On the other side we have
$$
\begin{array}{lll}
\delta(a)\delta(b)&=&\sum \eta_{h_{(1)}}(\tilde{g}_{(1)})f(S(h_{(2)})S(\tilde{g}_{(2)})\cdot \tilde{m})g_{(1)}\otimes \eta_{\tilde{h}_{(1)}}\otimes S(g_{(2)})\cdot m\otimes \tilde{h}_{(2)}\cdot \tilde{f}\\
\\
&=&\sum \psi(\tilde{g}_{(1)}h_{(1)})f(S(\tilde{g}_{(2)}h_{(2)})\cdot \tilde{m})g_{(1)}\otimes \eta_{\tilde{h}_{(1)}}\otimes S(g_{(2)})\cdot m\otimes \tilde{h}_{(2)}\cdot \tilde{f}\\
\\
&=&\sum f(S(\psi(\tilde{g}_{(1)}h_{(1)})\tilde{g}_{(2)}h_{(2)})\cdot \tilde{m})g_{(1)}\otimes \eta_{\tilde{h}_{(1)}}\otimes S(g_{(2)})\cdot m\otimes \tilde{h}_{(2)}\cdot \tilde{f}\\
\\
&=& \sum \psi(\tilde{g} h)f(\tilde{m}) g_{(1)}\otimes \eta_{\tilde{h}_{(1)}}\otimes S(g_{(2)})\cdot m \otimes \tilde{h}_{(2)}\cdot \tilde{f}.
\\
\\
&=& \delta{(ab)}
\end{array}
$$
It easy to check that the inverse of $\delta$ is
$$
\begin{array}{llll}
\delta^{-1}: & \hat{\mathcal{A}}({\mathcal{G}})\otimes \mathcal{A}(M^{\tau})& \to& \hat{\mathcal{A}}({\mathcal{G}})\otimes \mathcal{A}(M)\\ \\ & g\otimes \eta_{h} \otimes m \otimes  f &\mapsto &\sum g_{(1)}\otimes   \eta_{h_{(1)}} \otimes g_{(2)}\cdot  m\otimes S(h_{(2)})\cdot f
\end{array}
$$
There is an automorphism $\sigma$ of $\mathcal{G}$ such that $\psi(ab)=\psi(b\sigma(a))$, see \cite[Prop 3.12]{vd2}. Then
$$
\eta_{h}(x)=\psi(xh)=\psi(\sigma^{-1}(h)x)=\hat{\eta}_{\sigma^{-1}(h)}(x)
$$
and
$$
t\cdot\eta_{h}=t \cdot \hat{\eta}_{\sigma^{-1}(h)}= \hat{\eta}_{\sigma^{-1}(h)S(t)}= \eta_{h\sigma(S(t))}
$$
In the next step we use that
$$
\sum (\sigma(S(t)))_{(1)}\otimes(\sigma(S(t)))_{(2)}=\sum \sigma (S(t_{(2)}))\otimes S^{-1}(t_{(1)})
$$
see \cite[Prop 3.14]{vd2}.
Then we obtain
$$
\begin{array}{lll}
\displaystyle
\delta(t\cdot (g\otimes \eta_{h} \otimes m \otimes  f))&=& \delta(\displaystyle\sum t_{(1)} g\otimes t_{(4)}\cdot \eta_{h} \otimes t_{(2)}\cdot m \otimes t_{(3)}\cdot f)\\
\\
&=&\displaystyle\sum t_{(1)} g_{(1)}\otimes \eta_{(h\sigma(S(t_{(5)})))_{(1)}} \otimes S(t_{(2)}g_{(2)})t_{(3)}\cdot m \otimes(h\sigma(S(t_{(5)})))_{(2)} t_{(4)}\cdot f
\\
\\
&=&\displaystyle\sum t_{(1)} g_{(1)}\otimes \eta_{h_{(1)}\sigma(S(t_{(5)}))} \otimes S(g_{(2)})\epsilon(t_{(2)})\cdot m \otimes h_{(2)} S^{-1}(t_{(4)}))t_{(3)}\cdot f
\\
\\
&=&\displaystyle\sum t_{(1)} g_{(1)}\otimes \eta_{h_{(1)}\sigma(S(t_{(4)}))} \otimes S(g_{(2)})\cdot m \otimes h_{(2)} S^{-1}(t_{(3)}))t_{(2)}\cdot f
\\
\\
&=&\displaystyle\sum t_{(1)} g_{(1)}\otimes \eta_{h_{(1)}\sigma(S(t_{(2)}))} \otimes S(g_{(2)})\cdot m \otimes h_{(2)}\cdot f
\\
\\
&=&\displaystyle\sum t_{(1)} g_{(1)}\otimes t_{(2)}\cdot\eta_{h_{(1)}} \otimes S(g_{(2)})\cdot m \otimes h_{(2)}\cdot f
\\
\\
&=& t\cdot \delta(g\otimes \eta_{h} \otimes m \otimes  f)
\end{array}
$$
\end{proof}
\end{prop}

\begin{cor}
Let $M$ be a $\mathcal{G}$-module with good duality, then
$$
\hat{\mathcal{A}}(\mathcal{G})\otimes \hat{\mathcal{A}}(M) \simeq
\hat{\mathcal{A}}(\mathcal{G})\otimes \hat{\mathcal{A}}(M^{\tau})
$$
\end{cor}

\begin{rem}\label{mordif}
From Proposition \ref{modmodt} we take $\alpha$ and $\beta$. We know $\beta$ is a $\mathcal{G}$-module isomorphism, then by Proposition \ref{misoa} $\beta\otimes \alpha^{*}$
is a $\mathcal{G}$-module algebra isomorphism.
We have the following diagram
$$
\xymatrix{
\mathcal{A}(\mathcal{G}\otimes M)\ar[r]^{\beta\otimes \alpha^{*}}&
\mathcal{A}(\mathcal{G}\otimes M^{\tau})\\
\mathcal{A}(\mathcal{G})\otimes\mathcal{A}(M)\ar@{^{(}->}[u]&  \mathcal{A}(\mathcal{G})\otimes\mathcal{A}(M^{\tau})\ar@{^{(}->}[u]\\
\hat{\mathcal{A}}(\mathcal{G})\otimes\mathcal{A}(M)\ar@{^{(}->}[u]
&  \hat{\mathcal{A}}(\mathcal{G})\otimes\mathcal{A}(M^{\tau})\ar@{^{(}->}[u]
}
$$
In order to compare the map $\delta$ from Proposition \ref{isowow} with $\beta \otimes \alpha^{*}$ we take $h\in \mathcal{G}$ and $f\in M^*$ and compare $\alpha^{*}(\eta_{h}\otimes f)$ with $\sum \eta_{h_{(1)}}\otimes h_{(2)}\cdot f$ in $(\mathcal{G}\otimes M^{\tau})^{*}$
$$
\begin{array}{lll}
\alpha^{*}(\eta_{h}\otimes f)(k\otimes n) &=&(\eta_{h}\otimes f)(\alpha(k\otimes n))\\
\\
&=& (\eta_{h}\otimes f)(\displaystyle \sum k_{(1)}\otimes k_{(2)}\cdot n)) \\
\\
&=& \displaystyle \sum \psi (k_{(1)}h) f(k_{(2)}\cdot n)\\
\\
\displaystyle (\sum \eta_{h_{(1)}}\otimes h_{(2)}\cdot f)(k\otimes n) &=&\displaystyle \sum  \eta_{h_{(1)}}(k) (h_{(2)}\cdot f)(n)
\\
\\
&=&\displaystyle \sum \psi(kh_{(1)})f(S(h_{(2)})\cdot n)  \\
\end{array}
$$
\end{rem}

Let $A$ be a $\mathcal{G}$-module algebra and $M$ be a $\mathcal{G}$-module.
Define $\mathcal{A}(M)\otimes A$ with the following $\mathcal{G}$-module algebra structure:
$$
(m\otimes f \otimes a)(\tilde{m}\otimes \tilde{f} \otimes \tilde{a})=mf(\tilde{m})\otimes \tilde{f}\otimes a\tilde{a} \qquad \qquad t\cdot(m\otimes f \otimes a)=\sum t_{(1)}\cdot m\otimes t_{(3)}\cdot f \otimes t_{(2)}\cdot a
$$
If $M$ has good duality we can also consider $\hat{\mathcal{A}}(M)\otimes A$. As  $\hat{\mathcal{A}}(M)\otimes A$ is again a $\mathcal{G}$-module algebra and $\mathcal{G}$ is a $\mathcal{G}$-module with the regular action with good duality we consider
$$
 \hat{\mathcal{A}}(\mathcal{G})\otimes \hat{\mathcal{A}}(M)\otimes A
$$

\begin{cor}
Let $M$ be a $\mathcal{G}$-module with good duality, $A$ be a $\mathcal{G}$-module algebra then
$$
\delta_{A}=\delta\otimes \operatorname{id}_{A}:\hat{\mathcal{A}}(\mathcal{G})\otimes \hat{\mathcal{A}}(M) \otimes A \rightarrow
\hat{\mathcal{A}}(\mathcal{G})\otimes \hat{\mathcal{A}}(M^{\tau})\otimes A
$$
is an isomorphism of $\mathcal{G}$-module algebras.
\end{cor}

\subsection{Strong equivariant stabilization}
Let $(V_{1},W_{1},\phi_{1})$ and $(V_{2},W_{2},\phi_{2})$ be $(\mathcal{G},A)$-triples. Consider
$$
V=\left(\begin{array}{cc}0 & V_{1}\\ V_{2} & 0 \end{array}\right)  \quad W=\left(\begin{array}{cc}0 & W_{2}\\ W_{1} & 0 \end{array}\right) \quad V\otimes W =\left(\begin{array}{cc} V_{1}\otimes W_{1} & 0\\ 0 & V_{2}\otimes W_{2} \end{array}\right)
$$
and
$$
\phi: V\otimes W \to A \qquad \phi(v_{1}\otimes w_{1}+v_{2}\otimes w_{2})=  \phi_{1}(v_{1}\otimes w_{1})+ \phi_{2}(v_{2}\otimes w_{2})
$$
Let $\iota_{i}$ be the inclusion $i=1,2$,
\begin{equation}
\iota_{1}:W_{1}\underset{\tiny\phi_{1}}{\otimes}V_{1}\to W\underset{\tiny\phi}{\otimes}V\leftarrow   W_{2}\underset{\tiny\phi_{2}}{\otimes}V_{2}:\iota_{2}
\end{equation}
$$
 W\otimes_{A} V =\left(\begin{array}{cc} W_{2}\otimes_{A} V_{2} & 0\\ 0 & W_{1}\otimes_{A} V_{1} \end{array}\right)
$$
Let us check $\iota_{i}$ is an algebra morphism.
$$
 \iota_{1}(w\otimes v)=\left(\begin{array}{cc} 0 & 0\\ 0 & w \otimes v \end{array}\right) \qquad \iota_{2}(w\otimes v)=\left(\begin{array}{cc} w \otimes v & 0\\ 0 & 0  \end{array}\right)
$$
$$
\iota_{1}[(w\otimes v)(\tilde{w}\otimes \tilde{v})]=\iota_{1}(w\otimes \phi_{1}(v\otimes \tilde{w})\tilde{v})=\left(\begin{array}{cc} 0 & 0\\ 0 & (w \otimes v)(\tilde{w} \otimes \tilde{v}) \end{array}\right)
= \iota_{1}( w \otimes v) \iota_{1}(\tilde{w} \otimes \tilde{v})
$$

A functor $F:\mathcal{G}$-$\alg \to \mathcal{D}$ is {\emph{$(\mathcal{G},A)$-stable}} if $F(\iota_{1})$ and $F(\iota_{2})$ are isomorphisms for all $(\mathcal{G},A)$-triples $(V_{1},W_{1},\phi_{1})$ and $(V_{2},W_{2},\phi_{2})$.
We say that $F$ is {\emph{$\mathcal{G}$-stable}} if it is $(\mathcal{G},A)$-stable for all $A$.

\begin{defn}[Equivariant Morita equivalence]\label{eme}
Let $A$ and $B$ be unital $\mathcal{G}$-algebras. We say that $A$ is {\emph{equivariant Morita equivalent}} to $B$, denoted by $A\sim_{M}B$, if there exist left $\mathcal{G}$-modules $V$ and $W$ such that $V\in _{A}\mathcal{M}_{B}$ and $W\in _{B}\mathcal{M}_{A}$ together with $\mathcal{G}$-equivariant bimodule isomorphisms
$$
\beta: W\otimes _{A}V\to B \quad \quad \alpha: V\otimes_{B}W\to A.
$$
\end{defn}

\begin{rem}\label{reme}
There exisit $\gamma: W\otimes _{A}V\to B$  and  $\delta: V\otimes_{B}W\to A$ such that if we consider
$$
\varphi: V\otimes W \xrightarrow{\pi} V\otimes_{B}W\xrightarrow{\delta}  A \qquad \quad \psi:  W\otimes V \xrightarrow{\pi} W\otimes_{A}V\xrightarrow{\gamma} B
$$
then $\gamma$ is an isomorphism between $W\underset{\tiny\varphi}{\otimes}V$  and $B$ as $\mathcal{G}$-module algebras and  $\delta$  is an isomorphism between  $V\underset{\tiny\psi}{\otimes}W$ and $A$ as $\mathcal{G}$-module algebras.
\end{rem}

\begin{thm}
Let $F:\mathcal{G}$-$\alg\rightarrow B$ be a $\mathcal{G}$-stable functor. Let $A$, $B$ be $\mathcal{G}$-module algebras.
If $A\sim_{M}B$ then $F(A)\simeq F(B)$.
\end{thm}
\begin{proof}
Consider $A\oplus B$ as a $\mathcal{G}$-module algebra
$$
g\cdot(a+b)=g\cdot a + g \cdot b \qquad (a+b)(\tilde{a}+\tilde{b})=a\tilde{a}+b\tilde{b}
$$
As $A\sim_{M}B$ there exist left $\mathcal{G}$-modules $V$ and $W$ such that
$$V\in _{A}\mathcal{M}_{B} \qquad W\in _{B}\mathcal{M}_{A} \qquad  W\underset{\tiny\varphi}{\otimes}V \simeq B \qquad  V\underset{\tiny\psi}{\otimes}W\simeq A$$
where $\varphi$ and $\psi$ are as Remark \ref{reme}.
We consider $V\in _{A\oplus B}\mathcal{M}$, $W\in \mathcal{M}_{A\oplus B}$, $\hat{\varphi}$ and $\hat{\psi}$ as follows
$$\begin{array}{c}
(a+b)\cdot v=a\cdot v\qquad w\cdot(a+b)=w\cdot a\\
\\
 \hat{\varphi}:V\otimes W \to A\oplus B \quad  \hat{\varphi}(v\otimes w)=\varphi(v\otimes w)+0\\
\\
 \hat{\psi}:W\otimes V \to A\oplus B \quad  \hat{\psi}( w\otimes v)=0+\psi(w\otimes v)
\end{array}
$$
We obtain that $(V,W,\hat{\varphi})$ and $(W,V,\hat{\psi})$ are $(\mathcal{G},A\oplus B)$-triples. Observe that
$$
W\underset{\tiny\hat{\varphi}}{\otimes}V \simeq W\underset{\tiny\varphi}{\otimes}V \simeq B \qquad
\quad  V\underset{\hat{\tiny\psi}}{\otimes}W\simeq V\underset{\tiny\psi}{\otimes}W\simeq A
$$
We take
$$
\hat{V}=\hat{W}=\left(\begin{array}{cc}0 & V \\ W & 0 \end{array}\right)  \quad \hat{V}\otimes \hat{W} =\left(\begin{array}{cc} V \otimes W & 0\\ 0 & W\otimes V \end{array}\right)
$$
and
$$
\hat{\phi}(v\otimes w + \tilde{w}\otimes \tilde{v})=\hat{\varphi}(v\otimes w)+\hat{\psi}(\tilde{w}\otimes \tilde{v})
$$
As $F$ is a $\mathcal{G}$-stable functor $F(\iota_{1})$ and $F(\iota_{2})$ are isomorphisms where
$$
W \underset{\hat{\tiny\varphi}}{\otimes}V \xrightarrow{\iota_{1}} \hat{W}\underset{\hat{\tiny\phi}}{\otimes}\hat{V}\xleftarrow{\iota_{2}}   V\underset{\tiny\hat{\psi} }{\otimes} W
$$
Then $F(A)\simeq F(B)$.
\end{proof}

\subsection{Weak equivariant stabilization}
Let $(V,W,\phi)$ be a $\mathcal{G}$-triple and $A$ be a $\mathcal{G}$-module algebra.  Let
$$
\overline{V}=\left(\begin{array}{cc}0 & A \\ A\otimes V & 0 \end{array}\right)  \quad \overline{W}=\left(\begin{array}{cc}0 & W\otimes A\\ A  & 0 \end{array}\right) \quad \overline{V}\otimes \overline{W} =\left(\begin{array}{cc} A\otimes A & 0\\ 0 & A\otimes V \otimes W\otimes A  \end{array}\right)
$$
and
$$
\overline{\phi}: \overline{V}\otimes \overline{W }\to A \qquad \overline{\phi}(a\otimes b+c\otimes v\otimes w\otimes d)=  ab + \phi(v\otimes w)cd
$$
Observe such that $(\overline{V},\overline{W},\overline{\phi})$ is a
$(\mathcal{G},A)$-triple. Let $\overline{\iota}_{1}$ and
$\overline{\iota}_{2}$ be the inclusion,
$$
\overline{\iota}_{2}: A \to
\overline{W}\underset{\tiny\overline{\phi}}{\otimes}\overline{V}
\leftarrow  W\otimes A \otimes V:\overline{\iota}_{1}
$$

$$
 \overline{W}\otimes_{A} \overline{V} =\left(\begin{array}{cc} (W \otimes A)\otimes_{A} (A\otimes V)  & 0\\ 0 & A\otimes_{A} A
 \end{array}\right)=\left(\begin{array}{cc} W\otimes A \otimes V & 0\\ 0 & A
 \end{array}\right)
$$

Let $({V},{W},{\phi})$  be  a $\mathcal{G}$-triple. A functor
$F:\mathcal{G}$-$\alg \to \mathcal{D}$ is {\emph{weak stable with
respect to}} $({V},{W},{\phi})$ if $F(\overline{\iota}_{1})$ and
$F(\overline{\iota}_{2})$ are isomorphisms for all $\mathcal{G}$-module algebra $A$.

Let $M$ be a $\mathcal{G}$-module with good duality. Consider the
$\mathcal{G}$-triple $(\hat{M},M,\operatorname{ev})$.  A functor
$F:\mathcal{G}$-$\alg \to \mathcal{D}$ is {\emph{weak stable with
respect to $M$}} if it is weak stable with respect to
$(\hat{M}, M,\operatorname{ev})$. In other words, $F$ is weak
stable with respect to $M$ if the following morphisms of
$\mathcal{G}$-module algebras are isomorphisms through $F$

\begin{equation}\label{i1i2}
\iota_{2}: A \rightarrow \left(\begin{array}{cc} \hat{\mathcal{A}}(M)\otimes A & 0 \\
0 & A \end{array}\right) \leftarrow  \hat{\mathcal{A}}(M)\otimes A
:\iota_{1}
\end{equation}

A functor $F:\mathcal{G}$-$\alg \to \mathcal{D}$ is {\emph{weak $\mathcal{G}$-stable}} if it is weak stable with respect to $\mathcal{G}$ considered as a $\mathcal{G}$-module with the regular action.

\begin{ex}
Let $M=\C[\N]$ with the trivial action of $\mathcal{G}$. Consider $\chi_{i}:\N \rightarrow \C$ as $\chi_{i}(j)=1$ if $i=j$ and $\chi_{i}(j)=0$ if $i\neq j$.
Take
$$\hat{M}=\{ f\in M^{*}: \quad f =\displaystyle \sum_{i}a_{i}\chi_{i}:  \mbox{$a_{i}\in \C$ $a_{i}\neq 0$ for a finite amount of $i\in\N$}\} $$
then $M$ is a module with good duality and $\hat{\mathcal{A}}(M)=M_{\infty}(\C)$.
A functor $F:\mathcal{G}$-$\alg \to \mathcal{D}$ is weak stable with respect to $M$ if and only if $F$ is $M_{\infty}$-stable.
\end{ex}

\begin{ex} Let $\mathcal{X}$ be an infinity set.
We take $M=\C[\mathcal{X}]$ with the trivial action of $\mathcal{G}$ and
$$
\hat{M}=\{ f\in M^{*}: \quad f =\displaystyle \sum_{x}a_{x}\chi_{x}:  \mbox{$a_{x}\in \C$ $a_{x}\neq 0$ for a finite amount of $x\in\mathcal{X}$}\} $$
then $M$ is a module with good duality and $\hat{\mathcal{A}}(M)=M_{\mathcal{X}}(\C)$.
A functor $F:\mathcal{G}$-$\alg \to \mathcal{D}$ is weak stable with respect to $M$ if and only if $F$ is $M_{\mathcal{X}}$-stable.
\end{ex}

\begin{ex} Let $G$ be a countable group. Let $M=\C[G]$ with the regular action of $\mathcal{G}=\C G$. Take
$$\hat{M}=\{ f\in M^{*}: \quad f =\displaystyle \sum_{g}a_{g}\chi_{g}:  \mbox{$a_{g}\in \C$ $a_{g}\neq 0$ for a finite amount of $g\in G$}\} $$
then $M$ is a module with good duality and  $\hat{\mathcal{A}}(M)=M_{G}(\C)$.
A funtor $F:\mathcal{G}$-$\alg \to \mathcal{D}$ is weak stable with respect to $M$ if and only if $F$ is $G$-stable, \cite{eakk}.
\end{ex}

\begin{thm}\label{geneq}
Let $\mathcal{G}$ be an algebraic quantum group and $M$ be a
$\mathcal{G}$-module with good duality. Let $F:\mathcal{G}$-$\alg\to
\mathcal{D}$ be a weak stable functor with respect to $M^{\tau}$.
The functor
$$
\hat{F}:\mbox{$\mathcal{G}$-$\alg$} \to \mathcal{D} \qquad A\mapsto F(\hat{\mathcal{A}}({\mathcal{G}})\otimes A)
$$
is weak stable with respect to  $M$.
\end{thm}
\begin{proof}
We have to prove the morphisms \eqref{i1i2} are isomorphisms through $\hat{F}$. Observe $\hat{F}(\iota_{2})=F(\operatorname{id}_{\hat{\mathcal{A}}(\mathcal{G})}\otimes \iota_{2})$ and
$$
\xymatrix{\hat{\mathcal{A}}({\mathcal{G}})\otimes A \ar[rr]^-{\operatorname{id}_{\hat{\mathcal{A}}(\mathcal{G})}\otimes \iota_{2}} \ar[dd]&& \mbox{$\left(\begin{array}{cc}
\hat{\mathcal{A}}({\mathcal{G}}) \otimes \hat{\mathcal{A}}(M)\otimes A & 0 \\
0 & \hat{\mathcal{A}}({\mathcal{G}}) \otimes A \end{array}\right)$}
\ar[dd]^{\delta_{A}}\\
\\
\mbox{$\left(\begin{array}{cc}
\hat{\mathcal{A}}({\mathcal{G}})\otimes A  \otimes \hat{\mathcal{A}}(M^{\tau})& 0 \\
0 & \hat{\mathcal{A}}({\mathcal{G}}) \otimes A \end{array}\right)$}\ar[rr]_{\simeq}
 && \mbox{$\left(\begin{array}{cc}
\hat{\mathcal{A}}({\mathcal{G}}) \otimes \hat{\mathcal{A}}(M^{\tau})\otimes A & 0 \\
0 & \hat{\mathcal{A}}({\mathcal{G}}) \otimes A \end{array}\right)$}
}
$$
Idem for $\hat{F}(\iota_{1})$.
\end{proof}
\begin{rem}  If $G$ is a countable group then the assumption of Theorem \ref{geneq} for $\mathcal{G}=\C [G]$ is the same as \cite[Proposition 3.1.9]{eakk}.
\end{rem}

\section{Algebraic quantum kk-theory}\label{ekk}

 Let $\mathcal{X}$ be a set such that
$\operatorname{card}(\mathcal{X})=\operatorname{dim}_{\C}(\mathcal{G})\times \operatorname{card}(\N)$. Let $A$, $B$ be $\cG$-module algebras.
Define
$$
\kk^{\cG}(A,B)=\kk_{\mathcal{X}}(\hat{\mathcal{A}}(\cG)\otimes A, \hat{\mathcal{A}}(\cG)\otimes B)
$$
Consider the category $\mathfrak{KK}^{\cG}$ whose objects are the $\cG$-module algebras and the morphisms between $A$ and $B$ are the elements of $\kk^{\cG}(A,B)$. Recall that there exist an associative product
$$
\circ: \kk_{\mbox{\tiny$\cG$-$\alg$}}(B,C)\times \kk_{\mbox{\tiny$\cG$-$\alg$}}(A,B)\to \kk_{\mbox{\tiny$\cG$-$\alg$}}(A,C)
$$
which extends the composition of algebra homomorphisms. If $[\alpha]\in\kk_{\mbox{\tiny$\cG$-$\alg$}}(B,C)$ is an element represented by $\alpha: J^{n}(B)\to C^{\cS^n}$ and $[\beta]\in\kk_{\mbox{\tiny$\cG$-$\alg$}}(A,B)$ is an element represented by $\beta: J^{m}(A)\to B^{\cS^m}$ then $[\alpha]\circ[\beta]$ is an element of $\kk_{\mbox{\tiny$\cG$-$\alg$}}(A,C)$ represented by
\begin{equation}\label{lawc}
J^{n+m}(A)\xrightarrow{J^{n}(\beta)}J^{n}(B^{\cS^{m}})\xrightarrow{}J^{n}(B)^{\cS^{m}}\xrightarrow{\alpha^{\cS^{m}}} C^{\cS^{n+m}}.
\end{equation}
If we take $ \hat{\mathcal{A}}(\cG)\otimes A$, $ \hat{\mathcal{A}}(\cG)\otimes B$ and $  \hat{\mathcal{A}}(\cG) \otimes C$ instead of $A$, $B$ and $C$, we obtain the composition law in $\mathfrak{KK}^{\cG}$.

Let $j^{\cG}:\cG$-$\alg\r \mathfrak{KK}^{\cG}$ be the functor defined as the identity on objects and which sends each morphism of $\cG$-module algebras $f:A\to B$ to its class $[\id_{ \hat{\mathcal{A}}(\cG)}\otimes f]\in \kk^{\cG}(A,B)$.
The equivalence $\Omega: \mathfrak{KK}^{\cG}\to \mathfrak{KK}^{\cG}$ and distinguished triangles defined in \eqref{psec} give to $\mathfrak{KK}^{\cG}$ a triangulated category structure.

\begin{thm}\label{unikkh}
The functor $j^{\cG}:\mbox{$\cG$-$\alg$}\r \mathfrak{KK}^{\cG}$ is an
excisive, homotopy invariant, and weak $\cG$-stable
functor. Moreover, it is the universal functor for these properties. In other words,  if $\cT$
is a triangulated category and $R:\mbox{$\cG$-$\alg$} \r \cT$ together
with a class of morphisms $\{\overline{\partial}_{\operatorname{E}}:\operatorname{E}\in\cE\}$ is an
excisive, homotopy invariant and weak $\cG$-stable functor,
then there exists a unique triangle functor $\overline{R}:\mathfrak{KK}^{\cG}\r \cT$ such that the following diagram commutes
$$
\xymatrix{
 \mbox{$\cG$-$\alg$} \ar[r]^-{j^{\cG}}\ar[dr]_{R}& \mathfrak{KK}^{\cG}\ar@![d]^{\overline{R}}\\
&\cT
}
$$
\end{thm}
\proof This proof is similar to \cite[Thm 4.1.1]{eakk}. We use the same arguments taking $ \hat{\mathcal{A}}(\cG)$ instead of $M_{G}$.
Let $E: A\xrightarrow{f} B \xrightarrow{g} C$ be a weakly split extension in $\cG$-$\alg$. Define
$$
\begin{array}{ccc}
\partial^{\cG}_{\operatorname{E}}\in \hom_{\mathfrak{KK}^{\cG}}(\Omega C,
A)& = & \hom_{\mathfrak{KK}_{\mbox{\tiny{$\cG$-$\alg$}}}}( \hat{\mathcal{A}}(\cG)\otimes \Omega C,  \hat{\mathcal{A}}(\cG)\otimes A)\\
\\ &=& \hom_{\mathfrak{KK}_{\mbox{\tiny{$\cG$-$\alg$}}}}(\Omega( \hat{\mathcal{A}}(\cG)\otimes C),  \hat{\mathcal{A}}(\cG)\otimes A)
\end{array}
$$
as the morphism $\partial_{E'}$ associated to
$$
 \hat{\mathcal{A}}(\cG) \otimes A\r  \hat{\mathcal{A}}(\cG)\otimes B\r   \hat{\mathcal{A}}(\cG)\otimes C \quad \quad (\operatorname{E}')
$$
By Theorem \ref{xpuk} and Theorem \ref{geneq} the functor $j^{\cG}:\mbox{$\cG$-$\alg$}\to \mathfrak{KK}^{\cG}$ with the family $\{\partial^{\cG}_{\operatorname{E}}: \operatorname{E}\in\cE \}$ is an excisive, homotopy
invariant and weak $\cG$-stable functor.
Let $X:\cG$-$\alg \to \cT$ be another functor
which has the mentioned properties. By Theorem \ref{xpuk} there exists a unique triangle functor  $\overline{X}:\mathfrak{KK}_{\mathcal{X}}\r\cT$
such that the following diagram commutes
\begin{equation}\label{warm}
\xymatrix{
\mbox{$\cG$-$\alg$}\ar[ddr]_-{X}\ar[rr]^{j^{\cG}}\ar[dr]^-{j_{\mathcal{X}}}&& \mathfrak{KK}^{\cG}\ar@![ddl]^-{X'}\\
&\mathfrak{KK}_{\mathcal{X}}\ar[ur]\ar[d]_-{\overline{X}}&\\
&\cT&\\
}
\end{equation}
Define $X':\mathfrak{KK}^{\cG}\to \cT$ as $X'=X$ on objects and if $\alpha\in \kk^{\cG}(A,B)$ then
$$
X'(\alpha):=X(\iota_{B})^{-1}X(\iota'_{B})\overline{X}(\alpha)X(\iota'_{A})^{-1}X(\iota_{A}).
$$
Here $A\xrightarrow{\iota_{A}}\left(\begin{array}{cc} \hat{\mathcal{A}}(\mathcal{G})\otimes A & 0 \\
0 & A \end{array}\right)  \xleftarrow{\iota'_{A}}  \hat{\mathcal{A}}(\mathcal{G})\otimes A $ denotes the zig-zag between $A$ and $\hat{\mathcal{A}}(\mathcal{G})\otimes A$.
This is the unique possibility to define $X'$.
\qed

\section{Green-Julg Theorem for $\mathfrak{KK}^{\hp}$}\label{gjekk}

In this section we prove a version of the Green-Julg theorem for $\mathfrak{KK}^{\hp}$ when $\hp$ is a semisimple Hopf algebra.

Let $A$ be an $\cG$-module algebra. The {\emph{smash product algebra}}
$A\#\cG$ is $A\otimes \cG$ with the
following product
$$
(a\# h)(b \# k)=\sum a(h_{(1)}\cdot b) \# h_{(2)}k\quad \quad a,b \in A
\quad h,k\in \cG
$$
If $f:A\r B$ is a morphism of $\cG$-module algebras, we put
$$
f\#\cG: A\# \cG \r B \# \cG  \quad f\#\cG(a\# h)=f(a)\#h
$$
which is a morphism of algebras. Hence, we have a functor
$\#\ {\cG}:\cG$-$\alg \r \alg$.

\begin{prop}\label{isodes}
Let  $M$ be an $\cG$-module and $A$ an $\cG$-module algebra. The
following is an isomorphism of algebras
\begin{equation}\label{morisodes}
\phi: \mathcal{A}(M)\otimes (A \# \cG) \r (\mathcal{A}(M)\otimes
A)\#\cG \quad \phi(m\otimes f\otimes a\otimes h)=\sum
m\otimes (h_{1}\rightharpoondown f)\otimes a\# h_{2}.
\end{equation}
where $(h\rightharpoondown f)(x)=f(S(h)x)$
\begin{proof}
The inverse map is the following
\begin{equation}
\zeta:  (\mathcal{A}(M)\otimes A)\#\cG \r \mathcal{A}(M)\otimes (A \# \cG) \qquad
\zeta(m\otimes f\otimes a \otimes h)= \sum m\otimes (S^{-1}(h_{1})\rightharpoondown f) \otimes a \# h_{2}
\end{equation}
\end{proof}
\end{prop}

\begin{prop}\label{smakk}
There exists a unique functor ${\#}\ {\mathcal{G}}:\mathfrak{KK}^{\mathcal{G}}\r
\mathfrak{KK}$ such that the following diagram is commutative
$$
\xymatrix{
 \mbox{$\mathcal{G}$-$\alg$}\ar[r]^-{\# {\mathcal{G}}}\ar[d]_{j^{\hp}} &\alg \ar[d]^{j}\\
\mathfrak{KK}^{\mathcal{G}}\ar[r]_-{{\#} {\mathcal{G}}} &  \mathfrak{KK}
}
$$
\end{prop}
\proof
By Theorem \ref{unikkh} it is enough to prove $j(- \# {\mathcal{G}})$ is excisive,
homotopy invariant and weak $\mathcal{G}$-stable.
The properties are straightforward.
 \qed
\begin{prop}
The functor $\#\mathcal{G}:\mathcal{G}$-$\alg\to \alg$ is a triangle functor.
\end{prop}
\proof
 A distinguished triangle in $\mathfrak{KK}^{\mathcal{G}}$ is a diagram isomorphic to
$$
\Omega B\xrightarrow{j^{\mathcal{G}}(\iota)} P_{f}\xrightarrow{j^{\mathcal{G}}(\pi_{f})}
A\xrightarrow{j^{\mathcal{G}}(f)} B
$$
for some morphism of $\mathcal{G}$-algebras $f:A\r B$. That means, if it is isomophic to
\begin{equation}\label{htri}
\hat{\mathcal{A}}(\mathcal{G})\otimes\Omega B\xrightarrow{j_{\mbox{\tiny $\mathcal{G}$-$\alg$}}(\iota)} \hat{\mathcal{A}}(\mathcal{G})\otimes P_{f}\xrightarrow{j_{\tiny \mbox{$\mathcal{G}$-$\alg$}}(\pi_{f})}
\hat{\mathcal{A}}(\mathcal{G}) \otimes A\xrightarrow{j_{\mbox{\tiny $\mathcal{G}$-$\alg$}}(f)}\hat{\mathcal{A}}(\mathcal{G}) \otimes B
\end{equation}
in $\mathfrak{KK}_{\tiny \mbox{$\mathcal{G}$-$\alg$}}$.  The functor ${\#} \mathcal{G}:\mathfrak{KK}^{\mathcal{G}}\r \mathfrak{KK}$ sends the triangle \eqref{htri} to
\begin{equation}\label{ttrr}
(\hat{\mathcal{A}}(\mathcal{G})\otimes\Omega B)\#\mathcal{G} \to  (\hat{\mathcal{A}}(\mathcal{G}) \otimes P_{f})\# \mathcal{G} \to
(\hat{\mathcal{A}}(\mathcal{G})\otimes A)\# \mathcal{G} \to (\hat{\mathcal{A}}(\mathcal{G})\otimes B)\# \mathcal{G}
\end{equation}
which by Proposition \ref{isodes} and weak $\mathcal{G}$-stability is isomorphic to
\begin{equation}\label{seqt}
\Omega B\#\mathcal{G} \to  P_{f}\# \mathcal{G} \to A\# \mathcal{G} \to B \# \mathcal{G}
\end{equation}
As $P_{f}\# \mathcal{G}\simeq  P_{f\# \mathcal{G}}$ \eqref{seqt} is a distinguished triangle in $\mathfrak{KK}$.
\qed

Let $A$ be an $\mathcal{G}$-module algebra. Consider $\mathcal{G}\otimes A$ as a left
$\mathcal{G}$-module with the diagonal action. Also consider $\mathcal{G}\otimes A$ as
a right $A$-module with the regular action. In other words
$$
h\cdot(k\otimes a)= \sum h_{(1)}k \otimes h_{(2)}\cdot a \qquad (k\otimes
a)\cdot c = k\otimes ac
$$
It is easy to check that
\begin{equation}\label{dosac}
t\cdot((h\otimes a)\cdot c) = \sum (t_{(1)}\cdot (h\otimes a))\cdot(t_{(2)}\cdot
c)
\end{equation}
We define
$$
\en_{A}(\mathcal{G}\otimes A):= \{ f \in \en_{\C}(\mathcal{G}\otimes A) \mbox{
such that } f(k\otimes ac)= f(k\otimes a)\cdot c\}
$$
The $\mathcal{G}$-module structure in $\mathcal{G}\otimes A$ gives an $\mathcal{G}$-module algebra
structure in $\en_{\C}(\mathcal{G}\otimes A)$. It is easy to check that
$\en_{A}(\mathcal{G}\otimes A)$ is a sub-$\mathcal{G}$-module algebra of $\en_{\C}(\mathcal{G}\otimes A)$.

Consider $\hat{\mathcal{A}}(\mathcal{G})\otimes A$ as an $\mathcal{G}$-module algebra.
We have the following homomorphism of $\mathcal{G}$-module algebras
\begin{equation}\label{isoends}
T:\hat{\mathcal{A}}(\mathcal{G})\otimes A \r \en_{A}(\mathcal{G}\otimes A) \qquad \quad
T(f\otimes a)(h\otimes b)=f(h)\otimes ab
\end{equation}
If $\mathcal{G}=\hp$ is a finite dimensional Hopf algebra, then \eqref{isoends} is an isomorphism.

\begin{thm}\label{gjhalg}
Let $\hp$ be a semisimple Hopf algebra. The functor given by the trivial action
$ \tau: \mathfrak{KK} \r \mathfrak{KK}^{\hp}$ is left adjoint to the functor given by
the smash product $\# \hp :\mathfrak{KK}^{\hp}\r\mathfrak{KK}$.
In particular there is a natural isomorphism
$$
\kk^{\hp}(A^{\tau},B)\simeq \kk(A,B\# \hp) \quad \quad A\in \alg \quad
B\in \mbox{$\hp$-$\alg$}
$$
\end{thm}
\proof It is enough to prove that there exist
$\overline{\alpha}_{A}\in \kk(A, A^{\tau}\# \hp)$  and $\overline{\beta}_{B}\in
\kk^{\hp}(( B\# \hp)^{\tau},B)$
such that
 $$
 A^{\tau}\xrightarrow{\tau(\overline{\alpha}_{A})}(
 A^{\tau}\# \hp)^{\tau}\xrightarrow{\overline{\beta}_{A^{\tau}}}
 A^{\tau} \quad \mbox{ and } \quad
  B \#\hp \xrightarrow{\overline{\alpha}_{B\# \hp}}
 (B\# \hp)^{\tau}\#\hp \xrightarrow{\#(\overline{\beta}_{B})} B\# \hp
 $$
are the identities in  $\kk^{\hp}(A^{\tau},A^{\tau})$ and
$\kk(B\# \hp,B\# \hp)$ respectively.
As $\hp$ is semisimple there exists a right integral element $t\in \hp$
such that
$$
\epsilon(t)=1 \quad \quad th=\epsilon(h)t \quad \quad \forall h \in \hp
$$
We consider $\psi$ the right funtional such that $\psi(t)=1$,  the existence of invariants is guaranteed by \cite[Prop 5.1]{vd2}.
If we take the map $$R:\hp \rightarrow \hp \qquad R(z)=\sum\psi(t_{1}z)t_{2}$$
we can check $R*\id=\id =\id *R$ then $R$ is the antipode $S$. Then
\begin{equation}\label{gur}
S(z)=\sum\psi(t_{1}z)t_{2} \quad   z_1\otimes z_2=\sum S^{-1}(t_{3})\otimes \psi(t_{1}z)S^{-1}(t_{2})
\end{equation}

Define
$$
\alpha_{A}: A \r  A^{\tau}\# \hp = A\otimes \hp \qquad \quad
\alpha_{A}(a)=a\otimes t
$$
It is an algebra morphism because $t$ is idempotent.
Let
$$
\begin{array}{llll}
\beta_{B}: &(B \# \hp)^{\tau} &\r &\mathcal{A}(\hp)\otimes B \\
&b\# h &\mapsto & \displaystyle \sum S^{-1}(h_{2}t_{3})\otimes \hat{\eta}_{t_{1}}\otimes S^{-1}(h_{1}t_{2})\cdot b = \mho_{b\#h}
\end{array}
$$

One checks that $\beta_{B}$ is an equivariant algebra homomorphism.
Consider
\begin{equation}\label{alesid}
\begin{array}{ccccc}
 A^{\tau}&\xrightarrow{\tau({\alpha}_{A})}&(A^{\tau}\#\hp)^{\tau} &\xrightarrow{ {\beta}_{A^{\tau}}}&
 \mathcal{A}(\hp)\otimes A^{\tau}\\
a &\mapsto & a\otimes t &\mapsto &\mho_{a\otimes t}
\end{array}
\end{equation}
We can see $\mho_{a\otimes t}$ as a map from $\hp$ to $\hp\otimes A$, we denote by $\hat{t}$ to another copy of $t$ in order to differentiate both
$$
\begin{array}{cclr}
\mho_{a\otimes t}(z)&=& \displaystyle \sum S^{-1}(\hat{t}_{3}) S^{-1}(t_{2})\otimes \psi(\hat{t}_{1}z) S^{-1}(\hat{t}_{2}) S^{-1}(t_{1})\cdot a& \\
\\
&=& \displaystyle \sum z_{1} S^{-1}(t_{2})\otimes z_{2} S^{-1}(t_{1})\cdot a& \eqref{gur}\\
\\
&=& \displaystyle \sum z_{1} S^{-1}(t_{2})\otimes \epsilon(z_{2}) \epsilon(t_{1}) a& \mbox{$h\cdot a = \epsilon(h)a$}\\
\\
&=& \displaystyle \sum z S^{-1}(t)\otimes a&\\
\\
&=& \displaystyle \sum \epsilon(z) t\otimes a& \\
\\
&=& (t\otimes \epsilon \otimes a)(z)
\end{array}
$$
As $\hp$ is unimodular we have $t={S}^{-1}(t)$. Because $\mho_{a\otimes t}=t\otimes \epsilon\otimes a$

the map \eqref{alesid} represents the identity in
$\kk^{\hp}(A^{\tau},A^{\tau})$ in the sense of \cite[Rem. 2.6.2]{eakk}.

It remains to prove that the following morphism represents the identity in
 $\kk(B \# \hp,B\# \hp)$
\begin{equation}
\begin{array}{ccccccc}
 B\#\hp &\xrightarrow{{\alpha}_{B\# \hp}}&(B\#
 \hp)^{\tau}\otimes \hp &\xrightarrow{ {\# \beta}_{B}}&
  (\mathcal{A}(\hp)\otimes B)\# \hp & \xrightarrow{{\zeta}}& \mathcal{A}(\hp)\otimes (B\#\hp)\\
b\#h &\mapsto & b\# h\otimes t &\mapsto &\mho_{b\# h}\# t &\mapsto & \Lambda_{b\#h}
\end{array}
\end{equation}
The morphism $\zeta$ is defined in the Proposition \ref{isodes}.
The morphism defined in \eqref{isoends} is an isomorphism
$$
T: \mathcal{A}(\hp)\otimes (B\#\hp)\rightarrow \en_{B\#\hp}(\hp\otimes B\# \hp)
$$
where $T(h\otimes f\otimes b\# k)(x\otimes a\#y)=hf(x)\otimes (b\#k)(a\#y)$.

We shall
prove that the map
$b\# h \mapsto T( \Lambda_{b\# h})$
represents to the identity.
Define
$$
\begin{array}{c}
\delta: \hp \otimes B \# \hp\r \hp \otimes B \# \hp\quad \gamma:
\hp \otimes B \# \hp \r \hp \otimes B \# \hp
\\
\\
\delta(x\otimes a \# y)= \sum x_{1}\otimes
x_{2}\cdot a\#  x_{3} y \qquad
\gamma(x\otimes a \# y)=\sum x_{1}\otimes S(x_{3})\cdot a \# S(x_{2})y
\end{array}
$$
It is easy to check that they are mutually inverse.
By $M_{\infty}$-stability the following morphism represents the identity
$$
B \# \hp \r \mathcal{A}(\hp)\otimes (B\# \hp) \quad \quad b\#h\mapsto
t\otimes \epsilon \otimes b\#h
$$
Hence
$$
B \# \hp \r \en_{B\#\hp}(\hp\otimes B\# \hp) \quad \quad b\#h\mapsto
\delta\circ T (t\otimes \epsilon \otimes b\#h)\circ \gamma = T(\Lambda_{b\# h})
$$
also represents the identity, here we use \cite[Prop 5.1.2]{CT}.\qed

\begin{rem}
The adjointness between $\tau$ and $\#\hp$ fails to hold at the algebra level. A counterexample is given in \cite[Example 5.2.4]{eakk}.
\end{rem}

\begin{cor}
Let $\hp$ be a semisimple Hopf algebra. Let $A$ be an $\hp$-module algebra, then
$$
\kk^{\hp}(\mathbb{C}, A)\simeq \kk(\mathbb{C}, A\# \hp)\simeq \kh(A\#\hp).
$$
\end{cor}

\section{Baaj Skandalis duality}\label{bsd}

In \cite{eakk} we present a duality property in equivariant algebraic kk-theory. Given a group $G$ and a $G$-algebra $A$, we define a $G$-graduation on $A\rtimes G$.
If we take the crossed product for $G$-graded algebra, $A\rtimes G\hat{\rtimes} G$, we obtain a $G$-algebra isomorphic to $M_{G}\otimes A$. Every $G$-stable functor identify in the range
$A\rtimes G\hat{\rtimes} G$ with $A$.
We dualize de concept of $G$-action to $G$-graduation. If $G$ is a group, the existence of a group $\hat{G}$ such that satisfies
$$
A\rtimes G \rtimes \hat{G} \simeq A
$$
is not garantized.

In the context of algebraic quantum groups we have a good framework for group duality. If $\mathcal{G}$ is an algebraic quantum group there exists an algebraic quantum group $\hat{\mathcal{G}}$ which is the dual of $\mathcal{G}$. If $A$ is a $\mathcal{G}$-module algebra and we take the smash product $A\#\mathcal{G}$. The dual group $\hat{\mathcal{G}}$ acts on $A\#\mathcal{G}$
by the following way
$$
f \rightharpoonup (a\# g) =  a \# f \rightharpoonup g  \quad \mbox{where $ f \rightharpoonup g= \sum f(g_{(2)})g_{(1)} $ }
$$
We can take the smash product $A\#\mathcal{G}\# \hat{\mathcal{G}}$. From \cite{ahma} we obtain the following result:
\begin{thm}\cite{ahma} \label{dualt}
Let $A$ be a $\mathcal{G}$-module algebra then
$$
(A\#\mathcal{G})\#\hat{\mathcal{G}} \simeq \hat{\mathcal{A}}(\mathcal{G})\otimes A
$$
\end{thm}
We obtain a similar theorem of Baaj-Skandalis duality
\begin{thm}
The functors $\# \mathcal{G}:\mathfrak{KK}^{\mathcal{G}}\rightarrow \mathfrak{KK}^{\hat{\mathcal{G}}}$ and $\# \hat{\mathcal{G}}:\mathfrak{KK}^{\hat{\mathcal{G}}}\rightarrow \mathfrak{KK}^{{\mathcal{G}}}$ are equivalences and
$$
kk^{\mathcal{G}}(A,B)\simeq kk^{\hat{\mathcal{G}}}(A\#\mathcal{G}, B\#\mathcal{G})
$$
\end{thm}
\begin{proof}
Apply Theorem \ref{dualt} and $\mathcal{G}$-stabilization of $\mathfrak{KK}^{\mathcal{G}}$.
\end{proof}
\bibliographystyle{plain}

\begin{thebibliography}{10}

\bibitem{AK}
R. Akbarpour and M. Khalkhali.
\newblock{Hopf algebra equivariant cyclic homology and cyclic homology of crossed product algebras}.
\newblock{J. Reine Angew. Math. 559 (2003), 137--152}.


\bibitem{CT}
G.~Corti\~nas and A.~Thom.
\newblock {Bivariant algebraic $K$-theory}.
\newblock {\em Journal f{\"u}r die Reine und Angewandte Mathematik (Crelle's
  Journal)}, 610:267--280, 2007.

\bibitem{ahma}
Drabant, Bernhard; Van Daele, Alfons and Zhang, Yinhuo
\newblock{Actions of multiplier Hopf algebras.}
\newblock{ Comm. Algebra 27 (1999), no. 9, 4117�V4172.}

\bibitem{eakk}
E.~Ellis.
\newblock{Equivariant algebraic $kk$-theory and adjointness theorems}.
\newblock{\em J. Algebra},
\newblock{398 (2014), 200--226.}

\bibitem{Gar1}
G. Grarkusha
\newblock{Algebraic Kasparov K-theory I}
\newblock{Doc. Math. 19 (2014) 1207-1269.}

\bibitem{Gar2}
G. Garkusha
\newblock{Universal bivariant algebraic K-theories.}
\newblock{J. Homotopy Relat. Struct. 8 (2013), no. 1, 67-116.}

\bibitem{GoJ}
P.G.~Goerss and J.F.~Jardine.
\newblock{Simplicial homotopy theory}
\newblock{Progress in Mathematics (Boston, Mass.), 1999.}

\bibitem{ght}
E.~Guentner, N.~Higson, J~Trout.
\newblock{Equivariant $E$-theory for $C\sp *$-algebras.}
\newblock{Mem. Amer. Math. Soc. 148 (2000), no. 703, viii+86 pp.}


\bibitem{Mac} S. Mac Lane.
\newblock{Categories for the working mathematician.}
\newblock{\em 2nd ed., Graduate Texts in Mathematics.}
\newblock{ 5. New York, NY: Springer, 1998.}


\bibitem{Mkk}
R.~Meyer.
\newblock{Categorical aspects of bivariant K-theory}
\newblock{Corti\~nas, Guillermo (ed.) et al., $K$-theory and noncommutative
    geometry. Proceedings of the ICM 2006 satellite conference, Valladolid,
    Spain, August 31--September 6, 2006. Z\"urich: European Mathematical Society
    (EMS). Series of Congress Reports, 1-39, 2008}

\bibitem{RC}
E. Rodriguez.
\newblock{Bivariant algebraic K-theory categories and a spectrum for $G$-equivariant bivariant algebraic $K$-theory}
\newblock{PhD Thesis (2017). Universidad de Buenos Aires.}
\bibitem{sm}
S. Montgomery. {\it Hopf Algebras and Their Actions on Rings}. CBMS Regional Conference Series in Mathematics {\bf 82} (1993).

\bibitem{vd}
A. Van Daele.
\newblock{Multiplier Hopf Algebras}
\newblock{Trans. Amer. Math. Soc. 342 (1994), 917-932.}

\bibitem{vd2}
A. Van Daele.
\newblock{An algebraic framework for group duality}
\newblock{Adv. Math. 140 (1998), no. 2, 323�V366}

\bibitem{cv1}
C. Voigt.
\newblock{Cyclic cohomology and Baaj-Skandalis duality.}
\newblock{J. K-Theory 13 (2014), no. 1, 115--145}

\bibitem{cv2}
C. Voigt.
\newblock{Equivariant periodic cyclic homology.}
\newblock{ J. Inst. Math. Jussieu 6 (2007), no. 4, 689--763.}

\bibitem{kh}
C. Weibel. {\it Homotopy Algebraic $K$-theory}. Contemporary
Math. {\bf 83} (1989)
461--488.

\end{thebibliography}

\end{document}